\documentclass[preprint]{imsart}

\usepackage[round]{natbib}
\usepackage{amsthm,amsmath}
\RequirePackage[colorlinks,citecolor=blue,urlcolor=blue]{hyperref}

\RequirePackage{hypernat}

\startlocaldefs
\usepackage{t1enc}

\usepackage[latin2]{inputenc}
\usepackage{amssymb}
\usepackage{amsmath,amsthm}
\usepackage{color}
\usepackage{mathrsfs}
\usepackage{latexsym}
\usepackage{epsfig}
\usepackage[a4paper,left=4cm,right=3cm,top=3cm,bottom=3cm]{geometry}
\usepackage{subfigure}
\usepackage{hyperref}

\newcommand{\E}{\mathbb{E}}
\newcommand{\psw}{\ensuremath{{p_{\text{\rm sw},1}}}}
\newcommand{\switchratio}{\ensuremath{\frac{\psw(X^n)}{p_{B,0}(X^n)}}}
\newcommand{\proj}{\ensuremath{\Pi_0}}
\newcommand{\projb}{\ensuremath{\Pi'_0}}
\newcommand{\restje}{\ensuremath{\text{\tt small}_{\vec{b}}(n)}}
\newcommand{\mbn}{\ensuremath{\mu_1^{(n)}}}
\newcommand{\man}{\ensuremath{\mu_0^{(n)}}}
\newcommand{\sqd}[2]{\ensuremath{\| #1 - #2 \|^2_2}}
\newcommand{\dgen}[2]{\ensuremath{d_{\text{\tt gen}}(#1 \|#2)}}
\newcommand{\dgenb}{\ensuremath{d_{\text{\tt gen}}}}
\newcommand{\prob}{\mathbb{P}}
\newcommand{\1}{\textbf{1}}

\newcommand{\commentout}[1]{}
\newcommand{\reals}{\mathbb{R}}

\newcommand{\cM}{\mathcal{M}}
\newcommand{\M}{\mathcal{M}}
\newcommand{\cE}{\mathcal{E}}
\newcommand{\cF}{\mathcal{F}}
\newcommand{\cD}{\mathcal{D}}
\newcommand{\cA}{\mathcal{A}}
\newcommand{\cB}{\mathcal{B}}

\newcommand{\sq}{d_{SQ}}
\newcommand{\std}{d_{ST}}
\newcommand{\renyi}{d_{R}}
\newcommand{\hellinger}{d_{H^2}}

\newtheorem{theorem}{Theorem}
\newtheorem{lemma}{Lemma}
\newtheorem{corollary}{Corollary}
\newtheorem{proposition}{Proposition}
\theoremstyle{definition}
\newtheorem{definition}{Definition}
\newtheorem{example}{Example}

\numberwithin{equation}{section}

\graphicspath{{./Afbeeldingen/}}

\endlocaldefs

\begin{document}

\begin{frontmatter}

\title{Almost the Best of Three Worlds: Risk, Consistency and
  Optional Stopping for the Switch Criterion in Nested Model Selection}
\runtitle{Almost the Best of Three Worlds}


\author{\fnms{St\'ephanie} \snm{van der Pas}\ead[label=e1]{svdpas@math.leidenuniv.nl}\thanksref{t1}}
\and
\author{\fnms{Peter} \snm{Gr\"unwald}\ead[label=e2]{Peter.Grunwald@cwi.nl}\thanksref{t1,t2}}\\
\printead{e1}\\ \printead{e2}

\thankstext{t1}{Research supported by the Netherlands Organization for Scientific Research.}
\thankstext{t2}{This research was supported
by NWO VICI Project 639.073.04.}

\runauthor{van der Pas and Gr\"unwald}

\begin{abstract}
We study the switch distribution, introduced by Van Erven et al. (\citeyear{Erven2012}),  applied to model selection and subsequent estimation. While
  switching was known to be strongly consistent, here we show that it
  achieves minimax optimal parametric risk rates up to a $\log\log n$
  factor when comparing two nested exponential families, partially
  confirming a conjecture by \cite{Lauritzen2012} and
  \cite{Cavanaugh2012} that switching behaves asymptotically like the
  Hannan-Quinn criterion. Moreover, like Bayes factor model selection
  but unlike standard significance testing, when one of the models
  represents a simple hypothesis, the switch criterion defines a {\em
    robust\/} null hypothesis test, meaning that its Type-I error
  probability can be bounded irrespective of the stopping rule.
  Hence, switching is consistent, insensitive to optional
  stopping and almost minimax risk optimal, showing that, Yang's
  (\citeyear{Yang2005}) impossibility result notwithstanding, it is
  possible to `almost' combine the strengths of AIC and Bayes factor
  model selection.

\end{abstract}

\end{frontmatter}

\section{Introduction}
We
consider the following standard model selection problem, where we have
i.i.d. observations $X_1, \ldots, X_n$ and we wish to select between two
nested parametric models,
\begin{equation}\label{theproblem}
\mathcal{M}_0 = \left\{p_{\mu} \mid \mu \in M_0 \right\} \quad \text{ and } \quad \mathcal{M}_1 = \{p_{\mu} \mid  \mu \in M_1 \}.
\end{equation}
Here the $X_i$ are random vectors taking values in some set
${\cal X}$,  $M_1
\subseteq \reals^{m_1}$ for some $m_1 > 0$ and ${\cal M}_0 =
\{p_\mu : \mu \in M_0 \} \subset \cM_1$ represents an
$m_0$-dimensional submodel of $\cM_1$,
where $0 \leq m_0 < m_1$. We may thus denote $\M_0$ as the `simple' and $\M_1$ as the
`complex' model.  We will assume that $\mathcal{M}_1$ is an
exponential family, represented as a set of densities on  ${\cal
  X}$ with respect to some fixed underlying measure, so that $p_{\mu}$
represents the density of the observations, and we take it to be given in its mean-value parameterization. As the notation indicates, we require, without loss of generality, that the parameterizations of ${\cal
  M}_0$ and ${\cal M}_1$ coincide, that is $M_0 \subset M_1$ is itself a
set of $m_1$-dimensional vectors, the final $m_1 - m_0$ components of
which are fixed to known values. We restrict ourselves to the case in which both $M_1$ and
the restriction of $M_0$ to its first $m_0$ components are products of
open intervals.

Most model selection methods output not just a decision $\delta(X^n)
\in \{0,1\}$, but also an indication $r(X^n) \in {\mathbb R}$ of the
strength of evidence, such as a $p$-value or a Bayes factor. As a
result, such procedures can often be interpreted as methods for
hypothesis testing, where ${\cal M}_0$ represents the {\em null\/}
model and ${\cal M}_1$ the alternative; a very simple example of our
setting is when the $X_i$ consist of two components $X_i \equiv
(X_{i1}, X_{i2})$, which according to $\cM_1$ are independent
Gaussians whereas under $\cM_2$ they can have an arbitrary bivariate
Gaussian distribution and hence can be dependent.  Since we allow
$\cM_0$ to be a singleton, this setting also includes some very simple,
classical yet important settings such as testing whether a coin is
biased (${\cal M}_0$ is the fair coin model, ${\cal M}_1$ contains all
Bernoulli distributions).

We consider three desirable properties of model selection methods: (a)
optimal worst-case risk rate of post-model selection estimation (with
risk measured in terms of squared error loss, squared Hellinger
distance, R\'enyi or Kullback-Leibler divergence); (b) consistency,
and, (c), for procedures which also output a strength of evidence
$r(X^n)$, whether the validity of the evidence is insensitive to
optional stopping under the null model. We evaluate the recently
introduced model selection criterion $\delta_{\text{sw}}$ based on the
switch distribution (Van Erven et al., \citeyear{Erven2012}) on properties (a), (b) and (c).

The switch distribution, introduced by\footnote{Matlab code for
  implementing model selection, averaging and prediction by the switch
  distribution is available at {\tt http://www.blackwellpublishing.com/rss}. In general run times are comparable
  to those of the corresponding Bayesian methods.}  Van Erven et al.,
(\citeyear{Erven2007}), was originally designed to address the
\emph{catch-up phenomenon}, which occurs when the best predicting
model is not the same across sample sizes.  The switch distribution
can be interpreted as a modification of the Bayesian predictive
distribution. It also has an MDL interpretation: if one corrects
standard MDL approaches \citep{Grunwald2007} to take into account that
the best predicting method changes over time, one naturally arrives at
the switch distribution. \cite{LheritierC15} describe a successful
practical application for two-sample sequential testing, related to
the developments in this paper but in a nonparametric context. We
briefly give the definitions relevant to our setting in
Section~\ref{sec:switchgen}; for all further details we refer to Van
Erven et al. (\citeyear{Erven2012}) and \ref{app:realswitch} in the Appendix. 

When evaluating any model selection method, there is a well-known
tension between properties (a) and (b) above: the popular AIC method
\citep{Akaike1973} achieves the minimax optimal parametric rate of
order $1/n$ in the problem above, but is inconsistent; the same holds
for the many popular model selection methods that asymptotically tend
to behave like AIC, such as $k$-fold (for fixed $k$) and
leave-one-out-cross-validation, the bootstrap and Mallow's $C_p$ in
linear regression \citep{Efron1986, Shao1997, Stone1977}.  On the
other hand, BIC \citep{Schwarz1978} is consistent in the sense that
for large enough $n$, it will select the smallest model containing the
`true' $\mu$; but it misses the minimax parametric rate by a factor of
$\log n$.  The same holds for traditional Minimum Description Length
(MDL) approaches \citep{Grunwald2007} and Bayes factor model selection
(BFMS) \citep{Kass1995}, of which BIC is an approximation.  This might
lead one to wonder if there exists a single method that is optimal in
both respects. A key result by \cite{Yang2005} shows that this is
impossible: any consistent method misses the minimax optimal rate by a
factor $g(n)$ with $\lim_{n \rightarrow \infty} g(n) = \infty$.

In Section~\ref{sec:mainresult} we show that, Yang's result
notwithstanding, the switch distribution allows us to get very close
to satisfying property (a) and (b) at the same time, at least in the
problem defined above (Yang's result was shown in a nested linear
regression rather than our exponential family context, but it does
hold in our exponential family setting as well; see the discussion at the end of 
Section~\ref{sec:Yang}). We prove that in our setting, the switch model
selection criterion $\delta_{\text{sw}}$ (a) misses the minimax
optimal rate only by an exceedingly small $g_{\text{sw}}(n) \asymp
\log \log n$ factor (Theorem~\ref{thm:risk}). Property (b), strong
consistency, was already shown by Van Erven et al. (\citeyear{Erven2012}). The factor
$g_{\text{sw}}(n) \asymp \log \log n$ is an improvement over the extra
factor resulting from Bayes factor model selection, which has
$g_{\text{\sc bfms}}(n) \asymp \log n$. Indeed, as discussed in the
introduction of \cite{Erven2012}, the catch-up phenomenon that the
switch distribution addresses is intimately related to the
rate-suboptimality of Bayesian inference. Van
Erven et al. (\citeyear{Erven2012}) show that,
while model selection based on switching is consistent, sequential
prediction based on model averaging with the switching method achieves
minimax optimal {\em cumulative\/} risk rates in general parametric
and nonparametric settings, where the cumulative risk at sample size
$n$ is obtained by summing the standard, instantaneous risk from $1$
to $n$. In contrast, in nonparametric settings, standard Bayesian
model averaging typically has a cumulative risk rate that is larger by
a $\log n$ factor. Using the cumulative risk is natural in sequential
prediction settings, but Van
Erven et al. (\citeyear{Erven2012}) left open the question of
how switching would behave for the more standard, instantaneous
risk. In contrast to the cumulative setting, we cannot expect to
achieve the optimal rate here by Yang's (\citeyear{Yang2005}) result,
but it is interesting to see that switching gets so close.

We now turn to the third property, robustness to optional
stopping. While consistency in the sense above is an asymptotic and
even somewhat controversial notion (see Section~\ref{sec:discussion}),
there exists a nonasymptotic property closely related to consistency
that, while arguably more important in practice, has received
relatively little attention in the recent statistical literature.
This is property (c) above, insensitivity to optional stopping. In
statistics, the issue was thoroughly discussed, yet never completely
resolved, in the 1960s; nowadays, it is viewed as a highly desirable
feature of testing methods by, for example, psychologists; see
\citep{Wagenmakers2007,Sanborn2014}. In particular, it is often argued \citep{Wagenmakers2007}
that the fixed stopping rule required by the classical Neyman-Pearson
paradigm severely and unnecessarily restricts the application domain
of hypothesis testing,
invalidating much of the $p$-values reported in the psychological
literature. Approximately 55\%  of psychologists admitted in a survey to deciding whether to collect more data after looking at their results to see if they were significant \citep{John2012}. 
We analyze property (c) in terms of \emph{robust null
  hypothesis tests}, formally defined in Section~\ref{sec:robustsign}.
A method defines a robust null hypothesis test if (1) it outputs
evidence $r(X^n)$ that does not depend on the stopping rule used to
determine $n$, and (2) (some function of) $r(X^n)$ gives a bound on
the Type-I error that is valid no matter what this stopping rule is.
Standard (Neyman-Pearson) null hypothesis testing and tests derived
from AIC-type methods are not robust in this sense. For example, such
tests cannot be used if the stopping rule is simply unknown, as is
often the case when analyzing externally provided data --- but this is
just the tip of an iceberg of problems with nonrobust tests. For an
exhaustive review of such problems we refer to \cite{Wagenmakers2007}
who builds on, amongst others, \cite{Berger1988} and \cite{Pratt1962}.

Now, as first noted by \cite{Edwards1963}, in simple versus composite
testing (i.e. when $\cM_0$ is a singleton), the output of BFMS, the
Bayes factor, does provide a robust null hypothesis test.  This is one
of the main reasons why for example, in psychology, Bayesian testing
is becoming more and more popular \citep{Dienes2011, Andrews2012},
even among `frequentist' researchers \citep{Sanborn2014}.  Our third
result, in Section~\ref{sec:robustsign}, shows that the same holds for
the switch criterion: if $\cM_0$ is a singleton, so that the problem
\eqref{theproblem} reduces to a simple versus composite hypothesis
test, then the evidence $r(X^n)$ associated with the switching
criterion has the desired robustness property as well and thus in this
sense behaves like the Bayes factor method. The advantage, from a
frequentist point of view, of switching as compared to Bayes is then
that switching is a lot more sensitive: our risk rate results directly
imply that the Type II error $(1-\text{power})$ of the switch
criterion goes to $0$ as soon as, at sample size $n$, the distance
between the `true' distribution $\mu_1$ and the null model,
i.e. $\inf_{\mu \in M_0} \| \mu - \mu_1\|^2_2$ is of order $(\log \log
n)/n$; for Bayes factor testing, in order for the Type-II error to
reach $0$, this distance must be of order $(\log n)/n$ (this was informally recognized by  \cite{LheritierC15}, who reported substantially larger power of switching as compared to the Bayes factor method in a sequential two-sample testing setting). 

Thus, for singleton ${\cal M}_0$, switching gives us `almost the best
of three worlds': minimax rate optimality up to a $\log \log n$ factor
(in contrast to BFMS), consistency (in contrast to AIC-type methods)
and nonasymptotic insensitivity to optional stopping (in contrast to
standard Neyman-Pearson testing) in combination with a small Type-II
error. For composite ${\cal M}_0$, we show in
Section~\ref{sec:robustsign} that nonasymptotic robustness to optional
stopping still holds, albeit only in a much weaker sense --- thus
pointing towards an obvious goal for future work, discussed in
Section~\ref{sec:discussion}: can we modify the switch distribution so
as to get full optional stopping robustness also for composite ${\cal
  M}_0$?

\paragraph{Organization} 
This paper is organized as follows. The switch criterion is introduced
in Section~\ref{sec:switchgen}. In Section~\ref{sec:rateopt}, we
provide some preliminaries: we list the loss/risk functions for which
our result holds, describe the sets in which the truth is assumed to
lie, and discuss the tension between consistency and
rate-optimality. Suitable post-model-selection estimators to be used
in combination with the switch criterion are introduced in
Section~\ref{sec:main}, after which our main result on the worst-case
risk of the switch criterion is stated. We also go into the
relationship between the switch criterion and the Hannan-Quinn
criterion in that section. In Section~\ref{sec:robustsign} we define
robust null hypothesis tests, give some examples, and show that
testing by switching has the desired nonasymptotic robustness to
optional stopping; in constrast, AIC does not satisfy such a property
at all and the Hannan-Quinn criterion only satisfies an asymptotic
analogue. We also provide some simulations that illustrate our results.
Section~\ref{sec:discussion} provides some additional discussion and
ideas for future work. All proofs are given in the Appendix.

\paragraph{Notations and Conventions}
We use $x^n = x_1, \ldots, x_n$ to denote $n$ observations, each
taking values in a sample space $\mathcal{X}$.  For a set of
parameters $M$, $\mu \in M$, and $x \in \mathcal{X}$, $p_{\mu}(x)$
invariably denotes the density or mass function of $x$ under the
distribution $\mathbb{P}_{\mu}$ of random variable $X$, taking values
in $\mathcal{X}$. This is extended to $n$ outcomes by independence, so
that $p_{\mu}(x^n) := \prod_{i=1}^n p_{\mu}(x_i)$ and ${\mathbb
  P}_{\mu}(X^n \in A_n)$, abbreviated to ${\mathbb P}_{\mu}(A_n)$,
denotes the probability that $X^n \in A_n$ for $X^n = X_1, \ldots,
X_n$ i.i.d. $\sim {\mathbb P}_{\mu}$. Similarly, ${\mathbb E}_{\mu}$
denotes expectation under ${\mathbb P_{\mu}}$.  As is customary, we
write $a_n \asymp b_n$ to denote $0 < \lim_{n \to \infty} \inf a_n/b_n
\leq \lim_{n \to \infty} \sup a_n/b_n < \infty$.  For notational
simplicity we assume throughtout this paper that whenever we refer to
a sample size $n$, then $n \geq 3$ to ensure that $\log \log n$ is defined
and positive. 

Throughout the text, we refer to standard properties of exponential
families without always giving an explicit reference; all desired
properties can be found, in precise form, in
\citep{BarndorffNielsen1978} and, on a less formal level, in
\citep[Chapter 18,19]{Grunwald2007}.

\section{\label{sec:switchgen}Model Selection by Switching}
The \emph{switch distribution} (Van Erven et al., \citeyear{Erven2007,Erven2012}) is a
modification of the Bayesian predictive distribution, inspired by
Dawid's (\citeyear{Dawid84}) `prequential' approach to statistics and
the related {\em Minimum Description Length\/} (MDL) Principle \citep{BarronRY98,Grunwald2007}. The
corresponding \emph{switch criterion} can be thought of as Bayes
factor model selection with a prior on meta-models, where each
meta-model consists of a sequence of basic models and associated
starting times: until time $t_1$, follow model $k_1$, from time $t_1$
to $t_2$, follow model $k_2$, and so on. The fact that we only need to
select between two nested parametric models allows us to considerably
simplify the set-up of Van Erven et al. (\citeyear{Erven2012}), who dealt with countably
infinite sets of arbitrary models. 

It is convenient to directly introduce the switch criterion as a
modification of the Bayes factor model selection (BFMS). Assuming
equal prior $1/2$ on each of the models $\M_0$ and $\M_1$, BFMS
associates each model $\M_k$, $k \in \{0,1\}$, with a {\em marginal
  distribution\/} $p_{B,k}$ with
\begin{equation}\label{eq:bayes}
p_{B,k}(x^n) :=
\int_{\mu \in M_k} \omega_k(\mu) p_\mu(x^n)d\mu,
\end{equation} where $\omega_k$ is
a prior density on $M_k$. It then selects model $\M_1$ if and only if
$p_{B,1}(x^n) > p_{B,0}(x^n)$.

The basic idea behind MDL model selection is to generalize this in the
sense that each model $\M_k$ is associated with {\em some\/}
 `universal' distribution $p_{U,k}$; one then picks the $k$
for which $p_{U,k}(x^n)$ is largest. $p_{U,k}$ may be  set to the
Bayesian marginal distribution, but other choices may be preferable in
some situations. Switching is an instance of
this; in our simplified setting, it amounts to associating $\M_0$ with
a Bayes marginal distributon $p_{B,0}$ as before. $p_{U,1}$ however is
set to the {\em switch distribution\/} $\psw$. This distribution
corresponds to a switch between models $\M_0$ and $\M_1$ at some
sample point $s$, which is itself uncertain; before point $s$, the
data are modelled as coming from $\M_0$, using $p_{B,0}$; after point $s$, they are
modelled as coming from $\M_1$, using $p_{B,1}$. 
Formally, we denote the strategy that switches from the simple to the
complex model after $t$ observations by $\bar p_t$; $\psw$ is then
defined as the marginal distribution by averaging $\bar p_t$ over $t$,
with some probability mass function $\pi$ (analogous to a Bayesian prior) over $t \in \{1,2,\ldots\}$:
\begin{align*}
\bar p_t(x^n)  & = p_{B,0}(x^{t-1}) \cdot p_{B,1}(x_{t},\ldots, x_n
\mid x^{t-1}) \\
\psw(x^n) & = \sum_{t=1}^{\infty} \pi(t) \bar p_t(x^n),
\end{align*} 
where switching at $t = 1$ corresponds to predicting with $p_{B,1}$ at each data point, and switching at any $t > n$ to predicting with $p_{B,0}$. We remind the reader that even for i.i.d. models,
$p_{B,1}(x_{t},\ldots, x_n \mid x^{t-1})$ usually depends on $x^{t-1}$ ---
the Bayes predictive distribution learns from data.  The model
selection criterion $\delta_{\text{sw}}$ mapping sequences of
arbitrary length to $k \in \{0,1\}$ is then defined, for each $n$, as
follows:
\begin{equation}\label{eq:deltasworig}
\delta_{\text{sw}}(x^n) = \begin{cases}
0 & \text{if } {\displaystyle \frac{\psw(x^n)}{p_{\text{B, 0}}(x^n)}
\leq 1} \\
1 & \text{if } {\displaystyle \frac{\psw(x^n)}{p_\text{B, 0}(x^n)} > 1}
\end{cases}.
\end{equation}
When defining $\psw$ it is sufficient to consider switching times that are equal to a
power of two. Thus, we restrict attention to `priors' $\pi$ on
switching time with support on $2^0, 2^1, 2^2, \ldots $. 
For our subsequent results to hold, $\pi$ should be such that
$\pi(2^i)$ decays like $i^{-\kappa}$ for some $\kappa > 1$. An example of such
a prior with $\kappa = 2$ is $\pi(2^i) = 1/((i+1)(i+2))$, $\pi(j) = 0$
for any $j$ that is not a power of 2.

To prepare for Theorem~\ref{thm:risk}, we instantiate the switch
criterion to the problem \eqref{theproblem}. We define $p_{B,1}$
as any distribution of the form (\ref{eq:bayes}) where $\omega_1$ is a
continuous prior density on $M_1$ that is strictly positive on all
$\mu \in M_1$. To define $p_{B,0}$ we need to take a slight detour,
because we parameterized $\cM_0$ in terms of an $M_0$ that has a fixed
value on its final $m_1-m_0$ components: it is an $m_0$-dimensional
family with an $m_1$-dimensional parameterization, so one cannot
easily express a prior on $\cM_0$ as a density on $M_0$. To overcome
this, we distinguish between the case that $m_0=0$ and $m_0 > 0$. In
the former case $M_0$ has a single element $\nu$, and we define
$p_{B,0} = p_\nu$. In the latter case, we define 
$\projb: M_0 \rightarrow \reals^{m_0}$ as the projection of $\mu \in
M_0$ on its first $m_0$ components, and $\projb(M_0) := \{\projb(\mu):
\mu \in M_0 \}$.  For $\mu \in M_0$, we define $p_{\projb(\mu)}=
p_{\mu}$, and we  then let $\omega_0$ be a continuous strictly
positive prior density
on $\projb(M_0)$, and we define $p_{B,0}(x^n) :=
\int_{\mu' \in \projb(M_0)} \omega_0(\mu') p_{\mu'}(x^n)d\mu'$.

Two important remarks are in order: first, the fact that we associate
$\M_1$ with a distribution incorporating a `switch' from $\M_0$ to
$\M_1$ {\em does not mean\/} that we really believe that data were
sampled, until some point $t$, according to $\M_0$ and afterwards
according to $\M_1$. Rather, it is suggested by prequential and MDL
considerations, which suggest that one should pick the model that
performs best in sequentially predicting data; and if the data are
sampled from a distribution in $\M_1$ that is not in $\M_0$, but quite
close to it in KL divergence, then $p_{\text{B,1}}$ is suboptimal for
sequential prediction, and can be substantially outperformed by
$\psw$. This is explained at length by Van
Erven et al. (\citeyear{Erven2012}), and Figure 1
in that paper especially illustrates the point. The same paper also explains how one can
use dynamic programming to arrive at an implementation that has the
same computational efficiency as computation of the standard Bayes
model selection decision. 

Second, the criterion (\ref{eq:deltasworig}) as defined here is not
$100\%$ equivalent to the special case of the construction of
Van
Erven et al. (\citeyear{Erven2012}) specialized to two models, but rather a
simplification thereof. We do this purely for ease of explanation:
varying the exponent $\kappa$ in the prior
$\pi(2^i) \propto i^{-\kappa}$ defined above --- which is a free
parameter of the switch distribution --- has a stronger effect on the
switch criterion than switching between the two versions of the switch
distribution. This is explained in the Appendix, where
we also explain why all our results continue to hold if we were to
follow the original construction; conversely, the strong consistency
result for the construction of Van
Erven et al. (\citeyear{Erven2012}) trivially continues to
hold for the criterion (\ref{eq:deltasworig}) used in the present
paper.

\section{\label{sec:rateopt}Rate-optimality of post-model selection estimators}

This section contains some background to our main result, Theorem \ref{thm:risk}. In Section \ref{sec:loss}, we first list the loss functions for which our main result holds, and define the CINECSI sets in which the truth assumed to lie. We then discuss the minimax parametric risk for our model selection problem in Section \ref{sec:riskdiscussion}. This section ends with a discussion on the generality of the impossibility result of \cite{Yang2005} in Section \ref{sec:Yang}.

\subsection{\label{sec:loss}Loss functions and CINECSI sets}
Let $\cM = \{ p_{\mu} \mid \mu \in M\}$ be an exponential family given
in its mean-value parameterization with $M \subset \reals^m$ a product
of $m$ open, possibly but not necessarily unbounded intervals for some $m > 0$; see the Appendix
for a formal definition of exponential families and mean-value
parameterizations. Note that we do not require the family to be
`full'; for example, the Bernoulli model with success probability $\mu
\in M_1 = (0.2,0.4)$ counts as an exponential family in our (standard)
definition.

Suppose that we measure the quality of a density $p_{\mu'}$ as an
approximation to $p_{\mu}$ by a loss function $L: M
\times M \rightarrow \mathbb{R}$. The standard definition of the (instantaneous) {\em
  risk\/} of estimator $\breve{\mu}:
\bigcup_{i > 0} {\cal X}^i \rightarrow M$  at sample size $n$, as defined relative to loss $L$, is given
by its expected loss,
$$
R(\mu, \breve{\mu}, n) =   \E_\mu\left[L(\mu, \breve\mu(X^n)) \right],
$$
where  $\mathbb{E}_{\mu}$ denotes expectation over $X_1,
\ldots, X_n \ \text{i.i.d.}\ \sim \mathbb{P}_{\mu}$. Popular loss
functions are: 
\begin{enumerate} \item The \emph{squared error loss}: $\sq(\mu', \mu) =
\|\mu' - \mu\|_2^2$;
\item The {\em standardized squared error loss\/} which is a version
  of the \emph{squared Mahalonobis distance}, defined as 
\begin{equation}\label{eq:mald}
\std(\mu' \| \mu ) := (\mu - \mu')^T I(\mu') (\mu - \mu'), 
\end{equation}
where $T$ denotes transpose, $I(\cdot)$ is the Fisher information matrix, and we view $\mu$ and $\mu'$ as column vectors;
\item The R\'enyi divergence of order $1/2$, defined as 
$$\renyi(\mu',\mu) =  - 2 \log \E_{\mu'}\left[ \left( {p_{\mu}(X)
      }/{p_{\mu'}(X)} \right)^{\frac{1}{2}}\right];$$
\item The squared
  Hellinger distance $\hellinger(\mu',\mu) =  2 \left(1-\ \E_{\mu'}\left[ \left( {p_{\mu}(X)
      }/{p_{\mu'}(X)} \right)^{\frac{1}{2}}\right] \right)$;
  \item The KL (Kullback-Leibler) divergence $D(p_{\mu'} \| p_{\mu})$, henceforth abbreviated
  to $D(\mu' \| \mu)$.
\end{enumerate}
 We note
that there is a direct relationship between the R\'enyi divergence and squared Hellinger distance:
\begin{equation}\label{eq:rh}
\hellinger(\mu',\mu) = 2 \left( 1- e^{-
    \renyi(\mu',\mu)/2} \right).
\end{equation}
In fact, as we show below, these loss functions are all equivalent (equal up to
universal constants) on CINECSI sets. Such sets will play an important
role in the sequel. They are defined as follows:
\begin{definition}[CINECSI]
A \emph{CINECSI (Connected, Interior-Non-Empty-Compact-Subset-of-Interior) subset} of a set $M$  is a connected
subset of the interior of $M$ that is itself compact and has  nonempty
interior.
\end{definition}
The following proposition is proved in the Appendix.
\begin{proposition}\label{prop:withinconstants}
  Let $M$ be the mean-value parameter
  space of an exponential family as above, and let $M'$ be a CINECSI
  subset of $M$. Then there exist positive constants $c_1, c_2, \dots, c_6$  such that for all $\mu, \mu' \in M'$,
\begin{equation}\label{eq:firstdispb}
c_1 \| \mu' - \mu  \|^2_2 \leq c_2 \cdot \std(\mu' \| \mu)
\leq \hellinger(\mu',\mu) \leq \renyi(\mu', \mu)
\leq D(\mu' \| \mu) \leq c_3 \| \mu' - \mu  \|^2_2.
\end{equation}
and for all $\mu' \in M', \mu \in M$ (i.e. $\mu$ is now not restricted to lie in $M'$),
\begin{equation}\label{eq:seconddispb}
\hellinger(\mu',\mu) \leq c_4 \| \mu' - \mu  \|^2_2 \leq c_5 \cdot \std(\mu' \| \mu)
\leq c_6 \| \mu' - \mu  \|^2_2.
\end{equation}
\end{proposition}
CINECSI subsets are a variation on the INECCSI sets of
\citep{Grunwald2007}. Our main result, Theorem \ref{thm:risk}, holds
for all of the above  loss functions, and for general `sufficiently
efficient' estimators. While the equivalence of the losses above on
CINECSI sets is a great help in the proofs, we emphasize that we never
require these estimators to be restricted to CINECSI subsets of $M$
--- although, since we require $M$ to be open, every `true' $\mu \in
M$ will lie in {\em some\/} CINECSI subset $M'$ of $M$, a
statistician who employs the model $\cM$ cannot know what this
$M'$ is, so such a requirement would be unreasonably strong.

\subsection{\label{sec:riskdiscussion}Minimax parametric risk}

We say that a quantity $f_n$ converges at rate
$g_n$ if $f_n \asymp g_n$. We say that an estimator
$\breve\mu $ is  
\emph{minimax-rate optimal} relative to a model ${\cal M} = \{ p_{\mu} \mid \mu \in M\}$ {\em restricted to a subset $M'\subset M$\/} if 
\begin{equation*}
\sup_{\mu \in M'} R(\mu, \breve\mu, n)
\end{equation*}
converges at the same rate as 
\begin{equation}\label{eq:mm}
\inf_{\dot\mu} \sup_{\mu \in M'} R(\mu, \dot\mu, n),
\end{equation}
where $\dot\mu$ ranges over all estimators of $\mu$ at sample size $n$, that is,  all
measurable functions from ${\cal X}^n$ to $M$. 

For parametric models, the minimax risk (\ref{eq:mm}) is typically of
order $1/n$ when $R$ is defined relative to any of the loss measures
defined in Section \ref{sec:loss} and $M'$ is an arbitrary CINECSI
subset of $M$ (Van der Vaart, \citeyear{vandervaart1998}) --- models
for which this holds include e.g. most
location families and all curved exponential families, which include as a special case all standard exponential families. For this reason, from now on we refer to $1/n$ as
the {\em minimax parametric rate}.  Note that, crucially, the
restriction $\mu \in M'$ is imposed only on the data-generating
distribution, not on the estimators, and, since we will require models
with open parameter sets $M$ such that for every $\delta > 0$, there
is a CINECSI subset $M'_{\delta}$ of $M$ with
$\sup_{\mu \in M} \inf_{\mu' \in M'_{\delta}} \| \mu - \mu' ||^2_2 <
\delta$, every possible $\mu \in M$ will also lie in some CINECSI
subset $M'_{\delta}$ that `nearly' covers $M_{\delta}$. This makes the
restriction to CINECSI $M'$ in the definition above a mild one. Still,
it is necessary: at least for the squared error loss, for most
exponential families (the exception being the Gaussian location
family), we have
$\inf_{\dot\mu} \sup_{\mu \in M'_{\delta}} R(\mu, \dot\mu, n) =
C_{\delta}/n$ for some constant $C_{\delta} > 0$, but the smallest
constant for which this holds may grow arbitrarily large as
$\delta \rightarrow 0$, the reason being that the determinant of the
Fisher information may tend to $0$ or $\infty$ as
$\delta \rightarrow 0$.
  
Now consider a model selection criterion $\delta: \bigcup_{i > 0}
{\cal X}^i \rightarrow \{0,1,\ldots, K-1\}$ that selects, for given
data $x^n$ of arbitrary length $n$, one of a finite number $K$ of
parametric models $\cM_0, \ldots, \cM_{K-1}$ with respective parameter
sets $M_0, \ldots, M_{K-1}$. One way to evaluate the quality of
$\delta$ is to consider the risk attained after first selecting a
model and then estimating the parameter vector $\mu$ using an
estimator $\breve{\mu}_k$ associated with each model ${\cal M}_k$.
This {\em post-model selection estimator\/} \citep{Leeb2005} will be
denoted by $\breve\mu_{\breve k}(x^n)$, where $\breve k$ is the index
of the model selected by $\delta$. The risk of a model selection
criterion $\delta$ is thus $R(\mu, \delta, n) = \E_\mu\left[L(\mu,
  \breve\mu_{\breve k}(X^n)) \right]$, where $L$ is a given loss
function, and its worst-case risk relative to $\mu$ restricted to
$M'_k \subset M_k$ is given by
\begin{equation}\label{eq:mmd}
\sup_{\mu \in M'_k} R(\mu, \delta, n) = \sup_{\mu \in M'_k} \E_\mu\left[L(\mu, \breve\mu_{\breve k}(X^n)) \right].
\end{equation}
We are now ready to define what it means for a model selection criterion to achieve the minimax parametric rate.

\begin{definition}
A model selection criterion $\delta$ \emph{achieves the minimax parametric rate\/} if
there exist estimators $\breve{\mu}_k$, one for each ${\cal M}_k$
under consideration, such that, for {\em every\/} CINECSI subset
$M'_k$ of $M$:
\begin{equation*}
\sup_{\mu \in M'_k} R(\mu, \delta, n) \asymp 1/n. 
\end{equation*}
\end{definition}
Just as in the fixed-model case, the restriction $\mu \in
M'_k$ is imposed only on the data-generating distribution, not on the
estimators. 

\subsection{\label{sec:Yang}The result of \cite{Yang2005} transplanted
to our setting}
In this paper, as stated in the introduction, we further specialize
the setting above to problem \eqref{theproblem} where we select
between two nested exponential families, which we shall always assume
to be given in their mean-value parameterization. To be precise, the
`complex' model $\mathcal{M}_1$ contains distributions from an
exponential family parametrized by an $m_1$-dimensional mean vector
$\mu$, and the `simple' model $\mathcal{M}_0$ contains distributions
with the same parametrization, where the final $m_1 - m_0$ components
are fixed to values $\nu_{m_0+1}, \ldots, \nu_{ m_1}$.  We
introduce some notation to deal with the assumption that $M_1$ and its
restriction of $M_0$ to its first $m_0$ components are products of
open intervals. Formally, we require that $M_1$ and $M_0$ are of the
form
\begin{align}\label{eq:parameters} M_1 &=
  (\zeta_{1, 1}, \eta_{1,1}) \times \ldots \times
  (\zeta_{1, m_1},\eta_{1, m_1}) \nonumber \\
M_0 &=
  (\zeta_{0, 1}, \eta_{0,1}) \times \ldots \times
  (\zeta_{0, m_0},\eta_{0,m_0}) \times \{\nu_{m_0+1} \} \times \ldots \times
\{\nu_{ m_1}\}
\end{align}
where, for $j = 1, \ldots,
m_0$, we have $- \infty \leq \zeta_{1, j} \leq \zeta_{0, j} <
\eta_{0, j} \leq \eta_{1, j} \leq \infty$; and for $j = m_0 + 1,
\ldots, m_1$, we have $- \infty \leq \zeta_{1, j} < \nu_{ j} <
\eta_{1,j} \leq \infty$.  

For example, $\cM_1$ could contain all normal distributions with mean $\mu$ and variance $\sigma^2$, with mean value parameters $\mu_1 = \mu^2 + \sigma^2$ and $\mu_2 = \mu$, and $M_1 = (0, \infty) \times (-\infty, \infty)$, while $\cM_0$ could contain all normal distributions with mean zero and unknown variance $\sigma^2$, so $M_0 = (0, \infty) \times \{0\}$.

\cite{Yang2005} showed in a linear regression context that a model
selection criterion cannot both achieve the minimax optimal parametric
rate and be consistent; a practitioner is thus forced to choose
between a rate-optimal method such as AIC and a consistent method such
as BIC. Inequality \eqref{eq:risk} below provides some insight into why this
{\em AIC-BIC dilemma\/} can occur. A similar inequality appears in
Yang's paper for his linear regression context, but it is still valid
in our exponential family setting, and the derivation --- which we now
give --- is essentially equivalent. 

To state the inequality, we need to relate $\mu_1 \in M_1$ to a
component in $M_0$. For any given $\mu_1 = (\mu_{1,1}, \ldots, \mu_{1,
  m_1})^T \in M_1$, we will define
\begin{equation}
\label{eq:mu0b}
\proj(\mu_1):= (\mu_{1,1}, \ldots, {\mu}_{1,m_0}, \nu_{m_0+1}, \ldots,
\nu_{m_1})^T
\end{equation} 
to be the {\em projection\/} of $\mu_1$ on $M_0$. The difference between $\proj$ of \eqref{eq:mu0b} and $\projb$ in Section \ref{sec:switchgen} is that $\proj$ is a function from $\reals^{m_1}$ to $\reals^{m_1}$, whereas $\projb$ is a function from $\reals^{m_1}$ to $\reals^{m_0}$; $\proj(\mu_1)$ and $\projb(\mu_1)$ agree in the first $m_0$ components. Note that
$\proj(\mu_1)$ obviously minimizes, among all $\mu \in M_0$, the
squared Euclidean distance $\| \mu - \mu_1 \|^2_2$ to $p_{\mu_1}$;
somewhat less obviously it also minimizes, among $\mu \in M_0$, the KL
divergence $D(p_{\mu_1} \| p_{\mu})$ \citep[Chapter 19]{Grunwald2007};
we may thus think of it as the `best' approximation of the `true'
$\mu_1$ within $M_0$; we will usually abbreviate $\proj(\mu_1)$ to
$\mu_0$.  \commentout{ Moreover, $\mu_0$ satisfies
\begin{equation}\label{eq:multicramer}
  \E_{\mu_1}\left[ (\mu_0 - \breve\mu_0(X^n))^T I(\mu_1) (\mu_0 - \breve\mu_0(X^n) )\right]
  \leq
  \E_{\mu_1}\left[ (\mu_1 - \breve\mu_1(X^n))^T I(\mu_1) (\mu_1 - \breve\mu_1(X^n) )
  \right] =  \frac{1}{n}.\end{equation}
Here the inequality 
is a consequence of (\ref{eq:mu0}) and (\ref{eq:mu0b}).
The equality is derived using that for exponential families, $I(\mu)$, the Fisher information matrix for one outcome, is equal to the
inverse covariance matrix of the sufficient statistic $\phi$
\citep{BarndorffNielsen1978}; and since  we are assuming (\ref{eq:meanisml}) and i.i.d. data, it is also equal to $n^{-1}$ times the inverse covariance matrix of the ML estimator. The equality now follows by definition of the covariance matrix.}

Let $A_n$ be the event that the
complex model is selected at sample size $n$.  Since $\cM_1$ is an
exponential family, the MLE $\widehat\mu_1$ is unbiased and $\widehat\mu_0$ coincides
with $\widehat\mu_1$ in the first  $m_0$ components,  so that
$\E_{\mu_1}\left[\mu _0 - \widehat{\mu}_0(X^n)\right] = 0$, and hence we can rewrite, for any $\mu_1 \in M_1$,
the squared error risk as 
\begin{align}
R(\mu_1, \delta, n)
&=  \E_{\mu_1}\left[ \1_{A_n} \sqd{\mu_1}{\widehat\mu_1(X^n)} + 
\1_{A^c_n} \sqd{\mu_1}{\widehat\mu_0(X^n)}  \right] \nonumber \\
&= \E_{\mu_1}\left[ \1_{A_n} \sqd{\mu_1}{\widehat\mu_1(X^n)} + 
\1_{A^c_n} \sqd{\mu_0}{ \widehat\mu_0(X^n)}   + 1_{A^c_n} 
\sqd{\mu_1}{\mu_0} \right] \nonumber \\
&\leq \E_{\mu_1}\left[ \sqd{\mu_1}{\widehat\mu_1(X^n)} + 
\sqd{\mu_0}{ \widehat\mu_0(X^n)}  \right] + \prob(A^c_n) 
\sqd{\mu_1}{\mu_0} \nonumber \\
& \leq 2 R(\mu_1,\widehat{\mu}_1, n)  + \prob(A^c_n) \sqd{\mu_1}{\mu_0}.  \label{eq:risk}
\end{align}
The first part of the proof of our main result,
Theorem~\ref{thm:risk}, extends this decomposition to general
estimators and loss functions.

The first term on the right of (\ref{eq:risk}) is of order $1/n$.  The
second term depends on the `Type-II error', i.e. the probability of
selecting the simple model when it is not actually true. A low
worst-case risk is attained if this probability is small, even if the
true parameter is close to $\mu_0$. This does leave the possibility
for a risk-optimal model selection criterion to incorrectly select the
complex model with high probability. In other words, a risk-optimal
model selection method may not be consistent if the simple model is
correct. The theorem by \cite{Yang2005}, arguing from decomposition
\eqref{eq:risk}, essentially demonstrates that it cannot be. Due to
the general nature of (\ref{eq:risk}), it seems likely that his result
holds in much more general settings: a procedure attains a low
worst-case risk by selecting the complex model with high probability,
which is excellent if the complex model is indeed true, but leads to
inconsistency if the simple model is correct. Indeed, we have shown in
earlier work that the dilemma is not restricted to linear regression,
but occurs in our exponential family problem (\ref{theproblem}) as
well as long as $\cM_0 = \{ \nu \}$ is a singleton (see
(Van der Pas, \citeyear{Vanderpas2013}) for the proof, which is a simple adaptation of
Yang's proof that, we suspect, can be extended to nonsingleton $\cM_0$
as well). Hence, as the switch criterion is strongly consistent
(Van
Erven et al. (\citeyear{Erven2012})), we know that the worst-case risk rate of the switch
criterion cannot be of the order $1/n$ in general.

\section{\label{sec:main}Main result} 
We perform model selection by using the switch criterion, as specified in Section
\ref{sec:switchgen}. After the model selection, we estimate the
underlying parameter $\mu$. We discuss post-model selection estimators
suitable to our problem in Section \ref{sec:suffeff}. We are then
ready to present our main result, Theorem \ref{thm:risk} in Section
\ref{sec:mainresult}, stating that the worst-case risk for the switch
criterion under the loss functions listed in Section \ref{sec:loss}
attains the minimax parametric rate up to a $\log \log{n}$ factor.

\subsection{\label{sec:suffeff}Post-model selection: sufficiently efficient estimators}
Our goal is to determine the worst-case rate for the switch criterion
applied to two nested exponential families, which we combine with an
estimator as follows: if the simple model is selected, $\mu$ will be
estimated by an estimator $\breve{\mu}_0$ with range $M_0$. If the
complex model is selected, the estimate of $\mu$ will be provided by
another estimator $\breve{\mu}_1$ with range $M_1$.  Our result will
hold for all estimators $\breve{\mu}_0$ and $\breve{\mu}_1$ that are
{\em sufficiently efficient\/}: 
\begin{definition}[sufficiently efficient]
The estimators
$\{\breve{\mu}_k \rightarrow {M}_k \mid k \in \{0,1\} \}$ are
\emph{sufficiently efficient} with respect to a divergence measure $\dgen{\cdot}{\cdot}$  if (with $\mu_0 = \proj(\mu_1)$ as in (\ref{eq:mu0b})), for every CINECSI subset $M'_1$ of $M_1$, 
there exists a constant $C > 0$ such that for all $n$,
\begin{equation}\label{eq:suffeff}
\sup_{\mu_1 \in M'_1} \E_{\mu_1} [\dgen{\mu_0}{\breve\mu_0}] \leq C
\cdot \sup_{\mu_1 \in M'_1} \E_{\mu_1}[\dgen{\mu_1}{\breve\mu_1}] \leq
\frac{C}{n}.
\end{equation}
\end{definition}
Note that this is a stronger requirement than just rate-optimality: we
additionally require that, if the estimate $\breve\mu_0$ is used on
data sampled from $\mu_1 \in {M}_1$ (`misspecification'), then still
$\breve\mu_0$ converges to $\mu_0$, the best approximation of $\mu_1$
within $M_0$ at rate $O(1/n)$. In the Appendix we
provide a detailed discussion of sufficiently efficient estimators by
means of several examples. In a nutshell, it turns out that for
(standardized) squared error and Hellinger distance, the MLE is either
sufficiently efficient (e.g. for the Gaussian and gamma families), or
can be made sufficiently efficient by trivial modifications. For the
same losses, Bayes MAP estimates based on proper priors are
sufficiently efficient without modification for nearly all exponential
families. For R\'enyi and KL divergences, MLEs can sometimes be problematic but Bayes MAP estimators are still usually sufficiently efficient.

\subsection{\label{sec:mainresult}Main result: risk of the switch criterion}\

We now present our main result, which states that for the
exponential family problem under consideration, the worst-case
instantaneous risk rate of $\delta_\text{sw}$ is of order
$(\log\log{n})/n$.  Hence, the worst-case instantaneous risk of
$\delta_\text{sw}$ is very close to the lower bound of $1/n$, while
the criterion still maintains consistency. 

The theorem holds for any of the loss functions listed in Section
\ref{sec:loss}. We denote this by using the generic loss function
$\dgenb$, which can be one of the following loss functions: squared
error loss, standardized squared error los, KL divergence, R\'enyi
divergence of order 1/2, or squared Hellinger distance. Apart from the
sufficiently efficient condition on $\breve{\mu}_0$ and $\breve\mu_1$,
there are two minor conditions on the priors used in defining the
switch distribution: assumption 2 below rules out the use of improper
prior densities, but will hold for any other prior normally considered
for exponential family inference. Assumption 3 requires that the prior
probability of switching at time $t=2^i$ is strictly decreasing and
not exponentially small in $i$. Since these priors are user-defined
and not dependent on the underlying true distribution, these
conditions can easily be satisfied in practice.
\begin{theorem} \label{thm:risk} Let ${\cal M}_0 = \{ p_{\mu} \mid \mu \in
  M_0 \}$ and ${\cal M}_1 = \{ p_{\mu} \mid \mu \in M_1 \}$ be nested
  exponential families in their mean-value parameterization, where
  $M_0 \subseteq M_1$ are of the form (\ref{eq:parameters}). Assume:
  \begin{enumerate}
  \item  $\breve\mu_0$ and $\breve\mu_1$ are sufficiently efficient
  estimators relative to the  chosen loss $\dgenb$;
  \item $\delta_\text{sw}$ is constructed with $p_{B,0}$ and $p_{B,1}$ defined as
  in Section \ref{sec:switchgen}  with priors $\omega_k$ that admit a strictly positive, continuous density;
  \item and $\psw$ is defined relative to a prior $\pi$ with support
  on $\{0,1,2,4,8,\ldots\}$ and $\pi(2^i) \propto i^{-\kappa}$ for some $\kappa >1$. 
  \end{enumerate}

Then for every CINECSI subset $M_1'$ of $M_1$, we have: 
\begin{equation*} \sup_{\mu_1 \in M_1'} R(\mu_1,
  \delta_\text{sw}, n) = O\left( \frac{\log\log n}{n} \right),
\end{equation*}
for $R(\mu, \delta_\text{sw}, n)$ the risk at sample size $n$ defined
relative to the chosen loss $\dgenb$. 
\end{theorem}

\begin{example}{\bf [Our Setting vs. Yang's]\ }\label{ex:yang} \rm
  \cite{Yang2005} considers model selection between two nested linear
  regression models with fixed design, where the errors are Gaussian
  with fixed variance. The risk is measured as the in-model squared
  error risk (`in-model' means that the loss is measured conditional
  on a randomly chosen design point that already appeared in the
  training sample).  Within this context he shows that every model
  selection criterion that is (weakly) consistent cannot achieve the
  $1/n$ minimax rate.  The exponential family result above leads one
  to conjecture that the switch distribution achieves $O( (\log \log
  n)/n)$ risk in Yang's setting as well. We suspect that this is so,
  but actually showing this would require substantial additional
  work. Compared to our setting, Yang's setting is easier in some and
  harder in other respects: under the fixed-variance, fixed design
  regression model, the Fisher information is constant, making
  asymptotic results hold nonasymptotically, which would greatly
  facilitate our proofs (and obliterate any need to consider CINECSI
  sets or undefined MLE's). On the other hand, evaluating the risk
  conditional on a design point is not something that can be directly
  embedded in our proofs.
\end{example}
\begin{example}{\bf [Switching vs. Hannan-Quinn]\ }{\label{ex:HQ} \rm
In their comments on
Van Erven et al. (\citeyear{Erven2012}), \cite{Lauritzen2012} and
\cite{Cavanaugh2012} suggested a relationship between the switch model
selection criterion and the criterion due to \cite{Hannan1979}. For
the exponential family models under consideration,
the Hannan-Quinn criterion with parameter $c$, denoted as
$\text{HQ}$, selects the simple model, i.e. $\delta_{ \text{HQ}}
(x^n) = 0$, if
\begin{equation*}
- \log p_{\widehat{\mu}_0}(x^n) < - \log p_{\widehat{\mu}_1}(x^n) +c\log\log n,
\end{equation*}
and the complex model otherwise. In their paper, Hannan and Quinn show
that this criterion is strongly
consistent for $c > 1$.  

As shown by \cite{Barron1999}, under some regularity conditions,
penalized maximum likelihood criteria achieve worst-case quadratic
risk of the order of their penalty divided by $n$. One can show
(details omitted) that this is also the case in our specific setting
and hence, that the worst-case risk rate of HQ for our problem is of
order $(\log\log n)/n$. Our main result, Theorem~\ref{thm:risk},
shows that the same risk rate is achieved by the switch distribution,
thus partially confirming the conjecture of \cite{Lauritzen2012} and
\cite{Cavanaugh2012}: HQ achieves the same risk rate as the switch
distribution and, for the right choice of $c$, is also strongly
consistent. This suggests that the switch distribution and HQ, at
least for some specific value $c_0$, may behave asymptotically
indistinguishably. The earlier results of Van der Pas (\citeyear{Vanderpas2013})
suggest that this is indeed the case if $\cM_0$ is a singleton; if
$\cM_0$ has dimensionality larger than $0$, this appears to be a
difficult question which we will not attempt to resolve here --- in this
sense the conjecture of \cite{Lauritzen2012} and \cite{Cavanaugh2012}
has only been partially resolved.

Because HQ and $\delta_\text{sw}$ have been shown to be both strongly
consistent and achieve the same rates for this problem, one may wonder
whether one criterion is to be preferred over the other. For this
parametric problem, HQ has the advantage of being simpler to analyze
and implement. The criterion $\delta_\text{sw}$ can however, be used
to define a robust hypothesis test as in Section~\ref{sec:robustsign}
below. As we shall see there, HQ is insensitive to optional stopping
in an asymptotic sense only, whereas robust tests such as the switch
criterion are insensitive to optional stopping in a much stronger,
nonasymptotic sense. Except for the normal location model, for which
the asymptotics are precise, the HQ criterion cannot be easily adapted
to define such a robust, nonasymptotic test. Another advantage of
switching is that it can be combined with arbitrary priors and applied
much more generally, for example when the constituting models are
themselves nonparametric \citep{LheritierC15}, are so irregular that standard asymptotics
such as the law of the iterated logarithm are no longer valid, or are
represented by black-box predictors such that ML estimators and the
like cannot be calculated. In all of these cases the switch criterion
can still be defined and --- given the explanation in the introduction
of Van
Erven et al. (\citeyear{Erven2012}) --- one may still expect it to perform well.  
}
\end{example}

\section{\label{sec:robustsign}Robust null hypothesis tests }
Bayes factor model selection, the switch criterion, AIC, BIC, HQ and
most model selection methods used in practice are really based on
thresholding the output of a more informative \emph{model comparison
  method}. This is defined as a function from data of arbitrary size
to the nonnegative reals. Given data $x^n$, it outputs a number
$r(x^n)$  between 0 and $\infty$ that is a deterministic function of
the data $x^n$. Every model comparison method $r$ and threshold $t$
has an associated model selection method  $\delta_{r,t}$ that outputs
1 (corresponding to selecting model $\mathcal{M}_1$) if $r(x^n) \leq t
$, and 0 otherwise.  As explained below, such model comparison methods
can often be viewed as performing a null hypothesis  test with $\cM_0$ the
null hypothesis, $\cM_1$ the alternative hypothesis and $t$ akin to a
significance level.

\textsl{Example 1 (BFMS):} The output of the Bayes factor model  comparison method is the posterior odds ratio
$r_{\rm Bayes}(x^n) = \mathbb{P}(\mathcal{M}_0|x^n)/\mathbb{P}(\mathcal{M}_1|x^n)$.  The associated model selection method (BFMS) with threshold $t$  selects model $\mathcal{M}_1$ if and only $r_{Bayes}(x^n) \leq t$.

\textsl{Example 2 (AIC):}  Standard AIC selects model $\mathcal{M}_1$
if $\log (p_{\widehat{\mu}_1}(x^n) / p_{\widehat\mu_0}(x^n)) > m_1-m_0$.  We may however consider more conservative versions of AIC that only select $\mathcal{M}_1$ if
\begin{equation}\label{eq:aiccons}
{\log (p_{\widehat{\mu}_1}(x^n) / p_{\widehat{\mu}_0}(x^n)) } - (m_1- m_0) \geq - \log t.
\end{equation}
We may thus think of AIC as a model comparison method that outputs the left-hand side of (\ref{eq:aiccons}), and that becomes a model selection method when supplied with a particular $t$. 

Now classical Neyman-Pearson null hypothesis testing requires the {\em
  sampling plan}, or equivalently, the {\em stopping rule}, to be
determined in advance to ensure the validity of the subsequent
inference. In the important special case of (generalized) likelihood
ratio tests, this even means that the sample size $n$ has to be fixed
in advance. In practice, greater flexibility in choosing the sample
size $n$ is desirable (\cite{Wagenmakers2007} provides sophisticated
examples and discussion). Below, we discuss hypothesis tests that
allow such flexibility by virtue of the property that their Type
I-error probability remains bounded irrespective of the stopping rule
used. These \emph{robust} null hypothesis tests are defined below. As
will be shown, whenever the null hypothesis $\cM_0 = \{ p_{\mu_0} \}$
is `simple', i.e. a singleton (simple vs. composite testing), both
Bayes factor model selection (BFMS) and the switch distribution define
such robust null hypothesis tests, whereas AIC does not and HQ does so
only in an asymptotic sense. As we argue in Section~\ref{sec:power},
the advantage of switching over BFMS is then that, while both share
the robustness Type-I error property, switching has significantly
smaller Type-II error (larger power) than BFMS when the `truth' is
close to $\cM_0$, which is a direct consequence of it having a smaller
risk under the alternative $\cM_1$. To make this point concrete, and
to indicate what may happen if ${\cal M}_0$ is not a singleton, we
provide a simulation study in Section~\ref{sec:simulation}.

\subsection{Bayes Factors with singleton $\cM_0$ are Robust under
  Optional Stopping}
In many cases, for each $0 < \alpha < 1$  there is an associated
threshold $t(\alpha)$, which is a strictly increasing function of
$\alpha$, such that for every $t \leq t(\alpha)$ we have that
 $\delta_{r,t}$ becomes a null hypothesis significance test (NHST)
with type-I error probability bounded by  $\alpha$. In particular,
then $\delta_{r,t(\alpha)}$ is a standard NHST with type-I error
bounded by $\alpha$.  For example, for AIC with $M_0 =\{0\}$ and
$M_1 = \reals$ representing the normal family of distributions with
unit variance, we may select $t(\alpha) = \exp(-2/z_{\alpha/2}^2)$, where
$z_{\alpha/2}$ is the upper $(\alpha/2)$-quantile of the standard
normal distribution. This results in the generalized likelihood ratio
test at significance level $\alpha$.

We say that model comparison method $r$ defines a \emph{robust null
  hypothesis test} for null hypothesis $\cM_0$ if for all $\mu_0 \in M_0$,
all $0 \leq \alpha \leq 1$, 
\begin{equation} \label{eq:defrobustnull}
{\mathbb P}_{\mu_0}(\exists n: \delta_{r,t(\alpha) }  (X^n) =1 ) \leq \alpha . 
\end{equation}
Hence, a test that satisfies (\ref{eq:defrobustnull}) is a valid NHST
test at each fixed significance level $\alpha$, independently of the stopping
rule used. If a researcher can obtain a maximum of $n$ observations,
the probability of incorrectly selecting the complex model will remain
bounded away from one, regardless of the actual number of observations
made.

It is well-known that Bayes factor model selection provides a robust
null hypothesis test if we set $t(\alpha) = \alpha$,
as long as $\cM_0$ is a singleton. In other words, we may view the
output of BFMS as a `robust' variation of the $p$-value.  This was
already noted by \cite{Edwards1963} and interpreted as a frequentist
justification for BFMS; it also follows immediately from the following
result. \\ \

\begin{theorem}[Special Case of Eq. (2) of
  \cite{shafer2011}]\label{thm:optstopnull}
  Let $\cM_0,\cM_1, M_0$ and $M_1$ be as in Theorem~\ref{thm:risk}
  with common support ${\cal X} \subset \reals^d$ for some $d > 0$.
  Let $(X_1,X_2, \ldots)$ be an infinite sequence of random vectors
  all with support ${\cal X}$, and fix
  two distributions, $\bar{\mathbb P}_0$ and $\bar{\mathbb P}_1$ on
  ${\cal X}^{\infty}$ (so that under both $\bar{\mathbb P}_0$ and
  $\bar{\mathbb P}_1$, $(X_1, X_2, \ldots)$ constitutes a random
  process). Let, for each $n$, $\bar{p}^{(n)}_j$ represent the
  marginal density of $(X_1,\ldots, X_n)$ for the first $n$ outcomes
  under distribution $\bar{\mathbb P}_j$, relative to some product
  measure $\rho^n$ on $(\reals^d)^n$ (we assume $\bar{\mathbb P}_0$
  and $\bar{\mathbb P}_1$ to be such that these densities exist).  
Then for all $\alpha \geq 0$,
\begin{equation*}
\bar{\mathbb P}_0 \left(\exists  n: \frac{\bar{p}^{(n)}_{0}(X^n)}{\bar{p}^{(n)}_1(X^n)} \leq \alpha  \right) \leq \alpha.
\end{equation*}
\end{theorem}
We first apply this result for Bayes factor model selection, with model priors $\pi_0 =
\pi_1 = 1/2$
, so that $r_{\rm Bayes}(x^n) =
\mathbb{P}(\mathcal{M}_0|x^n)/\mathbb{P}(\mathcal{M}_1|x^n) =
p_{B,0}(x^n)/p_{B,1}(x^n)$. We immediately see:
\begin{corollary}\label{cor:sbayes}
  If $M_0 = \{\mu_0\}$ represents a singleton null model, then
  $p_{B,0} = p_{\mu_0}$ so that, applying
  Theorem~\ref{thm:optstopnull} with
  ${\mathbb P}_{B,0}= {\mathbb P}_{\mu_0}$, we see from
  (\ref{eq:defrobustnull}) that, if we set $t(\alpha)= \alpha$, then
  Bayes factor model selection constitutes a robust hypothesis
  test for null hypothesis $\cM_0$.
\end{corollary}
What happens if ${\cal M}_0$ is not singleton? Full robustness would
require that (\ref{eq:defrobustnull}) holds for all $\mu_0 \in
M_0$. The simulations below show that this will in general not be the
case for Bayes factor model selection. Yet, the same reasoning as used
in Corollary~\ref{cor:sbayes} implies that we still have some type of
robustness in a much weaker sense, which one might call ``robustness
in prior expectation'' relative to prior $\omega_0$ on $M_0$. Namely,
we have for all $0 \leq \alpha \leq 1$:
\begin{equation} \label{eq:defweakrobustnull}
{\mathbb P}_{B,0}(\exists n: \delta_{r,t(\alpha) }  (X^n) =1 ) \leq \alpha , 
\end{equation}
where ${\mathbb P}_{B,0}$ is the Bayes marginal distribution under
prior $\omega_0$.  In other words, if the beliefs of a Bayesian who
adopts prior $\omega_0$ on model ${\cal M}_0$ were accurate, then the
BFMS method would still give robust $p$-values, independently of the
stopping rule. While for a subjective Bayesian, such a weak form of
robustness might perhaps still be acceptable, we will stick to the
stronger definition instead, equating `robust hypothesis tests' with
tests satisfying (\ref{eq:defweakrobustnull}) uniformly for all $\mu_0
\in M_0$.

\paragraph{Remark} In practice we may very well
  be interested in a significance level $\alpha_n$ that is a fixed
  function of the sample size $n$, i.e., given data $X^n$, 
we choose ${\cal M}_1$ iff the output of the model comparison method is
  larger than $t(\alpha_n)$. Both Bayesian and switch-based model
  comparison may be used in this manner, and
  Theorem~\ref{thm:optstopnull} still holds with $\alpha$ replaced by
  $\alpha_n$; we focus on
  the fixed $\alpha$ case for simplicity only. 

\subsection{AIC is not, and HQ is only Asymptotically Robust}
The situation for AIC is quite different from that for BFMS and
switching: for every function $t: (0, 1) \to \mathbb{R}_{>0}$, we
have, even for every \emph{single} $0 < \alpha < 1$, that
$\delta_{AIC,t(\alpha)}$ is \emph{not} a robust null hypothesis test
for significance level $\alpha$. Hence AIC  cannot be transformed into
a robust test in this sense. This can immediately be seen when
comparing a $0$-dimensional (fixed mean $\mu_0$) with a 1-dimensional
Gaussian location family $\cM_1$ 
(extension to general multivariate exponential families is
straightforward but involves tedious manipulations with the Fisher
information). Evaluating the left hand side of (\ref{eq:aiccons})
yields that $\delta_{AIC,t(\alpha)}$ will select the complex model if

\begin{equation}\label{eq:AICLIL} \left| \sum_{i=1}^n \tilde X_i \right| \geq \frac{\sqrt{2n}}{t(\alpha)},
\end{equation}
where the $\tilde X_i$ are variables with mean 0 and variance 1 if
$\mathcal{M}_0$ is correct. Hence, as a consequence of the law of the iterated logarithm
(see for example Van der Vaart (\citeyear{vandervaart1998})), with probability one,
infinitely many $n$ exist such that the complex model will be favored,
even though it is incorrect.

It is instructive to compare this to the HQ criterion, which, in this example, using the same notation
as in (\ref{eq:AICLIL}), selects the complex model if
\begin{equation*}
\left| \sum_{i=1}^n \tilde X_i \right | \geq \sqrt{2cn\log\log n}.
\end{equation*}
If $c > 1$ (the case in which HQ is strongly consistent), then this
inequality will almost surely not hold for infinitely many $n$, as
again follows from the law of the iterated logarithm. The reasoning
can again be extended to other exponential families, and we find that
the HQ criterion with $c> 1$ is robust to optional stopping in the
crude, asymptotic sense that the probability that there exist
infinitely many sample sizes such that the simple model is incorrectly
rejected is zero. Yet HQ does not define a robust hypothesis test in
the sense above: to get the numerically precise Type I-error bound
(\ref{eq:defrobustnull}) we would need to define $t(\alpha)$ in a
model-and sample-size-dependent manner, which is quite complicated in all cases except
the Gaussian location families where the asymptotics hold
precisely. We note that the same type of asymptotic robustness holds
for the BIC criterion as well.

\subsection{Switching with singleton $\cM_0$ is Robust under
  Optional Stopping}
\label{sec:power}
The main insight of this section is simply that, just like BFMS,
switching can be used as a robust null hypothesis test as well, as
long as $\cM_0$ is a singleton: we can view the switch distribution as
a model comparison method that outputs odds ratio
$r_{\rm sw}(x^n) = p_{B,0}(x^n)/p_{\rm sw,1}(x^n)$. Until now, we used
it to select model 1 if $r_{\rm sw}(x^n) \leq 1$. If instead we fix a
significance level $\alpha$ and select model $1$ if
$r_{\rm sw}(x^n) \leq \alpha$, then we immediately see, by applying
Theorem~\ref{thm:optstopnull} in the same way as for the Bayes factor
case, that $r_{\rm sw}$ constitutes a  robust null hypothesis
test as long as $\cM_0$ is a singleton model (of course, if we select
${\cal M}_1$ as soon as $r_{\rm sw}$ outputs $t \leq \alpha$, then
$\alpha$ is merely an upper bound on the Type-I error; the actual
value might even be lower, as illustrated in the simulations
below). Similarly -- at least if the priors involved in the switch
criterion are chosen independently of the stopping rule --- just like
BFMS, the result $r_{\rm sw}(x^n)$ of model comparison by switching
does not depend on the `sampling intentions' of the analyst, thus
addressing the two most problematic issues with Neyman-Pearson
testing --- at least for singleton ${\cal M}_0$. Yet, from a frequentist perspective, switching is preferable
to BFMS, since it has substantially better power (type-II error)
properties. As could already be seen from Yang's decomposition
(\ref{eq:risk}), there is an intimate connection between Type-II error
and the risk rate achieved by any model comparison method. Formally,
we have the following result, a direct corollary of
Theorem~\ref{thm:upperbound} of the Appendix, which is
itself a major building block of our main result
Theorem~\ref{thm:risk} (plug in $\gamma = \alpha^{-1}$ into
(\ref{eq:start}) to get the corollary):
\begin{corollary}\label{cor:obayes}
 Using the same notations and under the same conditions as Theorem~\ref{thm:risk}, 
 for any $\alpha > 0$, there exist constants $C_1, C_2 > 0$ such that, for every CINECSI subset $M'_1$ of $M_1$, for
  every sequence $\mu_1^{(1)}, \mu_1^{(2)}, \ldots$ of elements of
  $M'_1$ with for all $n$, $\inf_{\mu_0 \in M_0} \| \mu_1^{(n)} - \mu_0 \|_2^2 \geq C_1 (\log \log n)/n$, we have
\begin{equation}\label{eq:startc}
\prob_{\mu_1^{(n)}}\left( r_{\rm sw}(x^n) \geq \alpha \right) \leq \frac{C_2}{\log n}.
\end{equation}
\end{corollary}
Hence, for any fixed significance level, the power of testing by
switching goes to $1$ as long as the data are sampled from a
distribution $\mu^{(n)}_1$ in $M_1$ that is farther away from $M_0$
than order $(\log \log n)/n$; for BFMS, the power only goes to $1$ if
$\mu^{(n)}$ is farther away than order $O((\log n)/n)$.

Corollary~\ref{cor:obayes} holds for general ${\cal M}_0$ including
composite ones. Yet robustness to optional stopping (and
hence `almost the best of three worlds') only holds if ${\cal M}_0$ is
a singleton; if ${\cal M}_0$ is composite, then --- using again
the same argument as for the Bayes factor case (see
Corollary~\ref{cor:sbayes} and directly below) --- we immediately see
from Theorem~\ref{thm:optstopnull} that the much weaker `prior expected robustness' property (\ref{eq:defweakrobustnull}) still holds. But, the simulations below show that full robustness does 
fail if $\mu_0$ is `atypical', i.e. if it resides far out in the tails
of the prior $\omega_0$.  A major question for
future work is now obviously whether there exist versions of the
switch criterion that give a truly robust null hypothesis test
even under a composite null hypothesis $\cM_0$. We return to this question in Section~\ref{sec:discussion}. 

\subsection{\label{sec:simulation}Simulation study}
\begin{figure}
\begin{center}
\includegraphics[width = 0.7\textwidth]{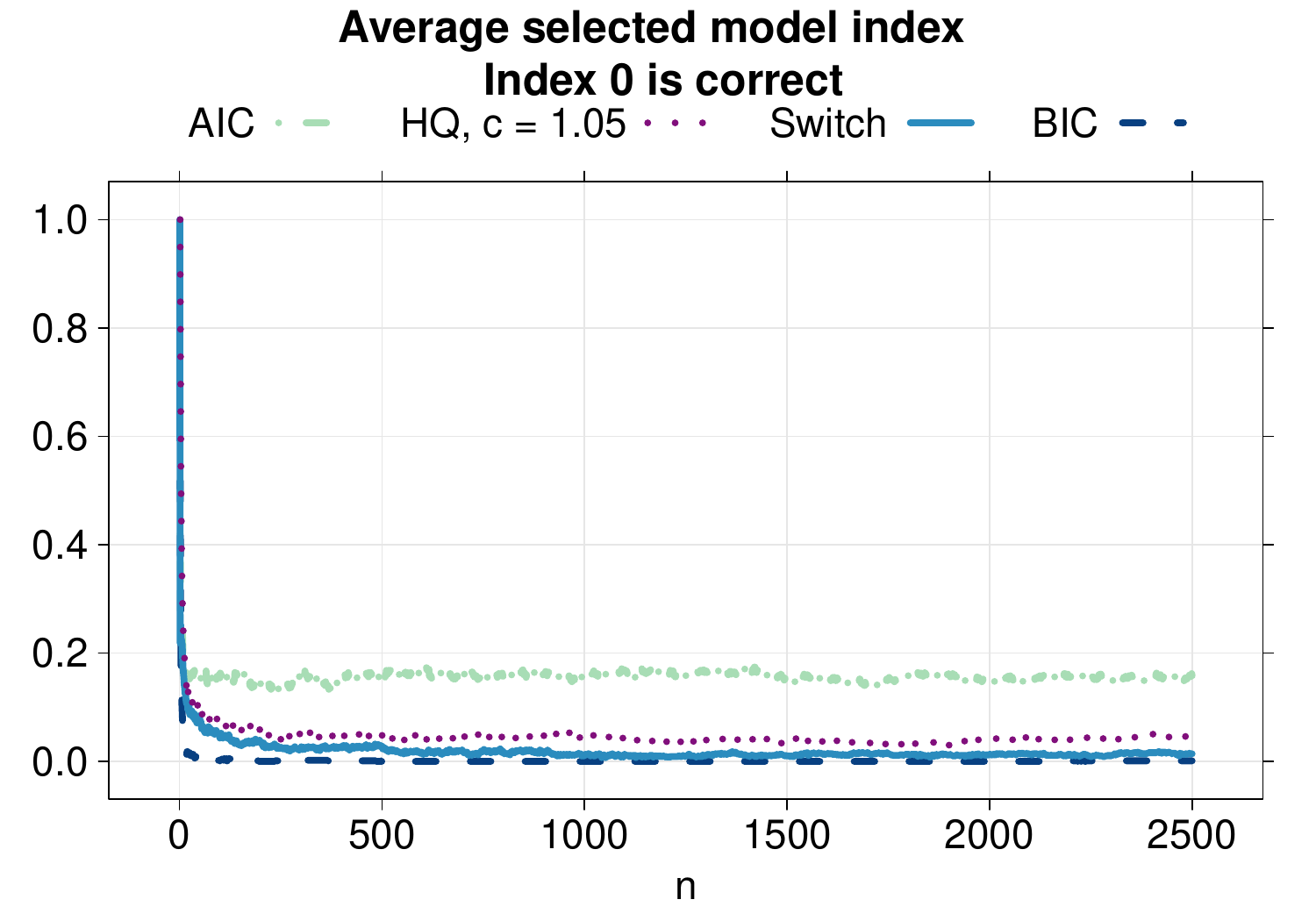}
\caption{$N = 1000$ data sets of length $n = 2500$ are generated from a standard normal distribution and the criteria are evaluated at each sample size. The figure shows the average selected model index (0 for $\mathcal{M}_1$, 1 for $\mathcal{M}_2$). The true index is 0. } 
\label{fig:consistency0}
\end{center}
\end{figure}

We now provide a simulation to illustrate the differences between AIC,
BIC, HQ and the switch criterion in terms of consistency, strong
consistency and robustness to optional stopping, illustrating the insights of the previous subsections. In each setting, two
of the following three models are compared:
\begin{itemize}
\item $\mathcal{M}_0 = \{\mathcal{N}(0,1)\}$.
\item $\mathcal{M}_1 = \{\mathcal{N}(\mu,1), \mu \in \mathbb{R}\}$, with a normal prior with mean zero and variance equal to 100 on $\mu$.
\item $\mathcal{M}_2 = \{\mathcal{N}(\mu, \sigma^2), \mu \in \mathbb{R}, \sigma \in \mathbb{R}_{>0}\}$, with a normal-inverse-gamma prior: $\mu | \sigma^2 \sim \mathcal{N}(0, C\times\sigma^2), \sigma^2 \sim IG(\alpha, \beta)$, with $C = 100, \alpha = 1, \beta = 1$.
\end{itemize}
To illustrate standard consistency, $\mathcal{M}_1$ and
$\mathcal{M}_2$ are considered. In the first setting, $\mathcal{M}_1$
is true. $N = 1000$ data sets of length $n = 2500$ are generated from
a standard normal distribution, and AIC, BIC, HQ with $c = 1.05$ and
$\delta_\text{sw}$ are evaluated at each sample size. The average
selected model index (0 for $\mathcal{M}_1$, 1 for $\mathcal{M}_2$) is
given in Figure \ref{fig:consistency0}.

\begin{figure}
\begin{center}
\includegraphics[width = 0.7\textwidth]{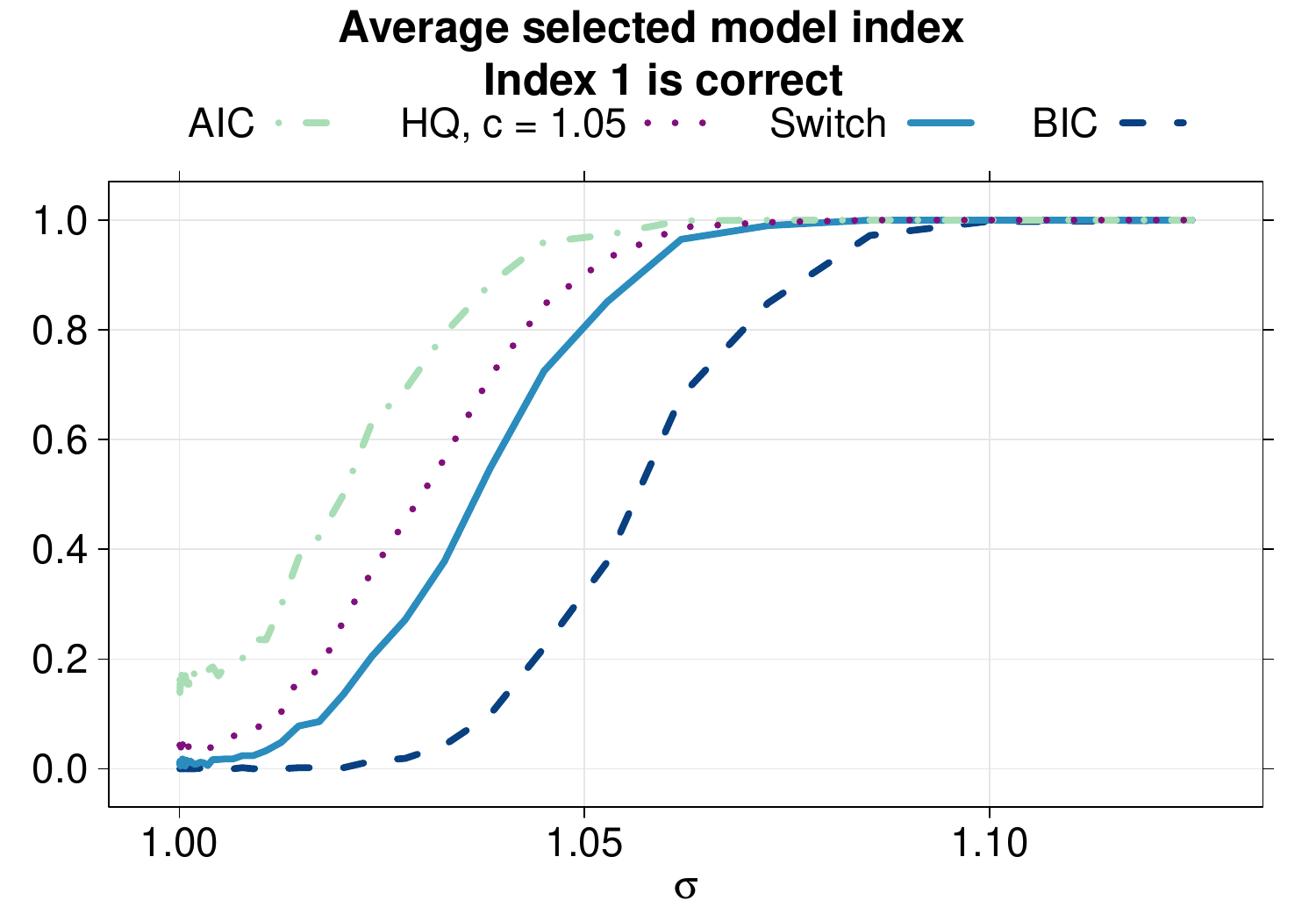}
\caption{$N = 1000$ data sets of length $n = 2500$ are generated from a normal distribution with mean 0 and variance $\sigma^2$ for a range of values of $\sigma$. The criteria are evaluated at $n = 2500$. The figure shows the average selected model index (0 for $\mathcal{M}_1$, 1 for $\mathcal{M}_2$). The true index is 1.  } 
\label{fig:consistency1}
\end{center}
\end{figure}

In the second setting, $\mathcal{M}_2$ is true. The data is generated from a normal distribution with mean 0 and a variance that is varied. For each value of $\sigma$, $N = 1000$ datasets of length $n = 2500$ are generated, and the four model selection criteria are evaluated at that sample size. The average selected model index  is given in Figure \ref{fig:consistency1}. 

The results are as expected. When the complex model is true, AIC is most likely to select it, at the cost of inconsistency when the simple model is true. BIC is the slowest to correctly select the complex model  and the first to correctly select the simple model. HQ and $\delta_\text{sw}$ show intermediate behaviour, HQ being slightly more likely to select the complex model.

\begin{figure}
\begin{center}
\includegraphics[width = 0.7\textwidth]{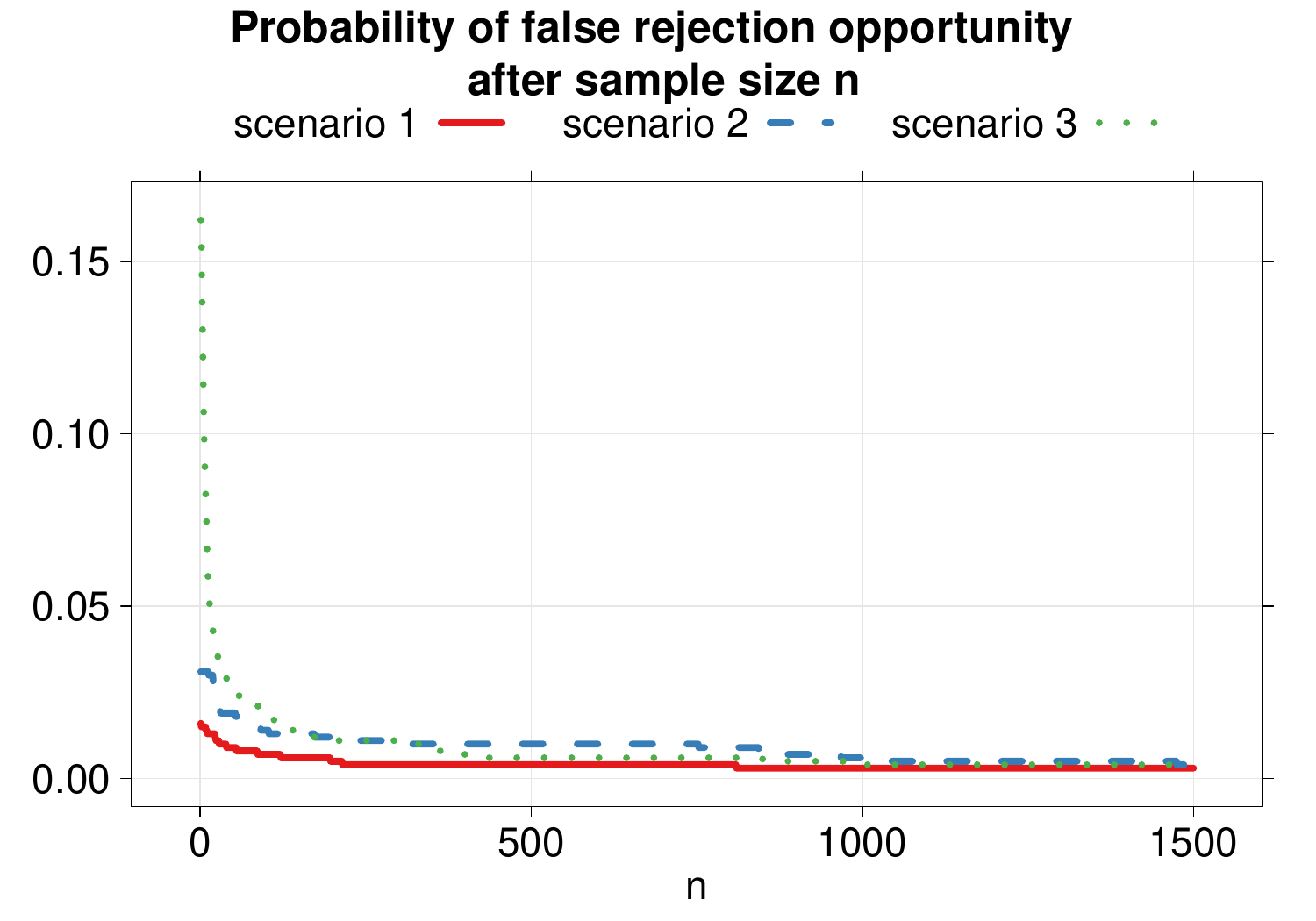}
\caption{ $N = 1000$ data sets of length $n_\text{max} = 10000$ in each scenario, from the simple model. The complex model is selected when $\delta_\text{sw}(x^n) > 20$. Estimated probability that there exists a model index after $n$ at which the complex model will be selected. Results shown up to $n = 1500$ for clarity. After $n = 1500$, the three curves are indistinguishable and all very close to zero. } 
\label{fig:optstopna}
\end{center}
\end{figure}

\begin{figure}
\begin{center}
\includegraphics[width = 0.7\textwidth]{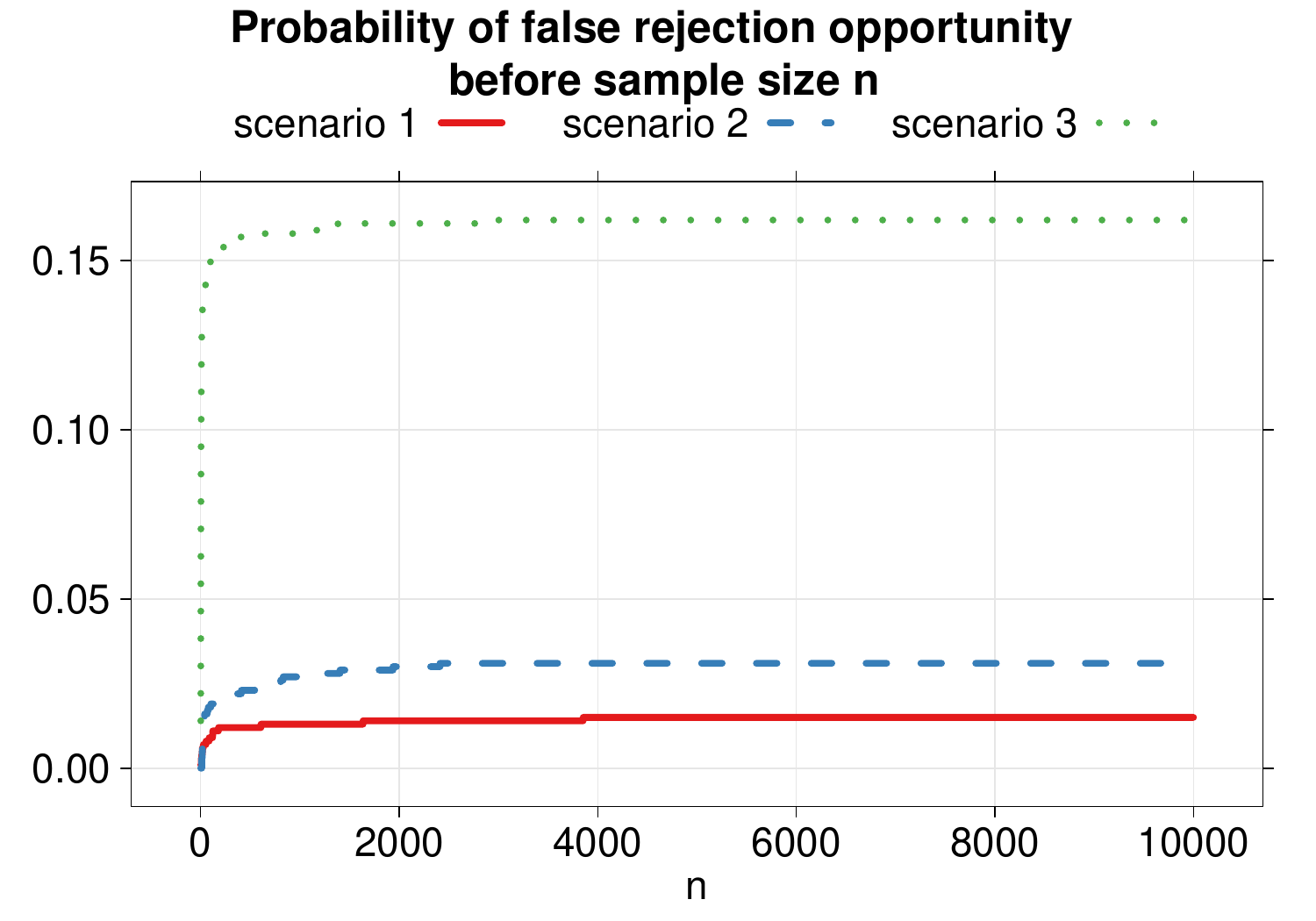}
\caption{  Setting as Figure \ref{fig:optstopna}. Estimated probability that there exists a model index before $n$ at which the complex model would have been selected. } 
\label{fig:optstopvoor}
\end{center}
\end{figure}

To illustrate strong consistency and optional stopping, three scenarios are considered:
\begin{enumerate}
\item $\mathcal{M}_0$ vs $\mathcal{M}_1$, data from a standard normal distribution (``scenario 1" --- Theorem \ref{thm:optstopnull}/Corollary~\ref{cor:sbayes} implies that switching defines a test that is  robust with respect to optional stopping).
\item $\mathcal{M}_1$ vs $\mathcal{M}_2$, data from a standard normal distribution (``scenario 2'', Theorem \ref{thm:optstopnull} does not only imply robustness, because null model is composite). 
\item $\mathcal{M}_1$ vs $\mathcal{M}_2$, data from a normal
  distribution with mean 35 and variance 1 (``scenario 3", Theorem
  \ref{thm:optstopnull} again does not imply robustness).
\end{enumerate}
We create $N = 1000$ data sets of length $n_\text{max} = 10000$ in
each scenario. We select the complex model when $\delta_\text{sw}$ is
larger than 20 (in terms of the robust $p$-value interpretation of
Theorem~\ref{thm:optstopnull}, this corresponds to a significance
level of $0.05$). We estimate two probabilities at each sample size
$n$:
\begin{itemize}
\item The probability that there will ever be a model index after $n$ at which the complex model will be selected (Figure \ref{fig:optstopna}), approximated by checking whether the complex model is selected at any sample size between $n$ and $3n_\text{max}$. 
\item The probability that there exists a model index before $n$ at which the complex model would have been selected (Figure \ref{fig:optstopvoor}).
\end{itemize}
Figure \ref{fig:optstopna} can be interpreted as a check whether
strong consistency holds --- if it does, then the probabilities should
converge to $0$ as $n \rightarrow \infty$. Van Erven et al.'s (2007)
theorem implies that strong consistency holds in all three scenarios,
and the graphs confirm this --- even though for scenario 3, in which
data comes from a $\mu \in M_0$ that is `atypical' under the prior, it
takes a bit longer --- illustrating that strong consistency is not a
uniform notion. The graph also illustrates that strong consistency can
be viewed as an asymptotic, nonuniform version of robustness to
optional stopping --- it implies that from some sample size (which may
be very large though) onwards, one will never again falsely reject no
matter how long one keeps sampling.

Figure~\ref{fig:optstopvoor} refers to nonasymptotic optional stopping:
in scenario 1, the conditions from Theorem \ref{thm:optstopnull} hold,
and indeed the figure shows that the probability that the complex
model is {\em ever\/} incorrectly selected even when optional stopping
is used, is bounded by 0.05 (the observed bound is 0.015). In
scenarios 2 and 3, the conditions from Theorem \ref{thm:optstopnull}
do not hold. In scenario 2, the behaviour of the switch criterion is
similar to scenario 1. However, in scenario 3, the probability of a
false rejection opportunity before sample size $n$ is not bounded by
0.05, but quickly goes to 0.15. We clearly see
that $\delta_\text{sw}$ is not  robust to optional
stopping in scenario 3.

When the simplest model is not a singleton, the choice of prior on the
model parameters (in scenarios 2 and 3 on $\mu$ in $\mathcal{M}_1$ and
on $(\mu, \sigma^2)$ in $\mathcal{M}_2$) affects the results. In both
scenario 2 and 3, $\delta_{\text{sw}}$ must still satisfy the weak,
prior-expected version of robustness (\ref{eq:defweakrobustnull}), as
we have seen in Section~\ref{sec:power}. In scenario 2, the prior is
centered at the data-generating value of zero and we do observe actual
robustness. In scenario 3 however, the prior is centered at zero while
the data is generated with a mean of 35, 3.5 standard deviations away
from the prior mean --- thus $\mu$ is `atypical' under the prior, and,
as the figure shows, nonasymptotic robustness is violated.

\section{Discussion and Future Work}
\label{sec:discussion}
In this paper we showed that switching combines near-rate optimality,
consistency and, for singleton ${\cal M}_0$, robustness to optional stopping.  We end the
paper by highlighting three issues which, we feel, need additional
discussion: first, the desirability of consistency; second, whether
there is anything `special' to the switch criterion as opposed to
other possible trade-offs between risk optimality and consistency; and
third, the limitations of switching in its current form.

\paragraph{Consistency}
Since the desirability of consistency, in the sense of finding the
smallest model containing the true distribution, is somewhat
controversial, let us discuss it a bit further.  The main argument
against consistency is made by those adhering to Box's maxim
`Essentially, all models are wrong, but some are useful'
\citep{Box1987}. According to some, the goal of model selection should
therefore not be to select a non-existing `true' model, but to obtain
the best predictive inference or best inference about a parameter
\citep{Burnham2004, Forster2000}. Another issue with consistency is
that it is a `nonuniform' notion, which in our context means that ---
as is indeed easy to see --- it is impossible to give a bound on the
probability under $\mathbb{P}_{\mu}$ of selecting the wrong model at
sample size $n$ that converges to $0$ uniformly for all $\mu \in
M$. This nonuniformity implies that consistency is of little practical
consequence for post-model selection inference \citep{Leeb2005}.

As to the first argument, one can reply that there do exist situations
in which a model can be correct, for example in the field of
extrasensory perception \citep{Bem2011}, in which it seems exceedingly
likely that the null model (expressing that no such thing exists) is
correct; another example is genetic linkage \citep{Gusella1983,
  Tsui1985}. The second argument is more convincing, but only to argue
that even if consistency holds, a method may not be very useful in
practice. It does not contradict that consistency can sometimes be a
highly desirable (but never the only highly desirable) property --- we
feel that this is the case whenever we are not purely interested in
prediction but instead are also seeking to find out whether a certain
structural relationship (e.g. dependence between variables) holds or
not.

Going one step further, it seems a good idea to study model selection
methods not in terms of the asymptotic, nonuniform notion of
consistency but instead by a more tangible finite-sample analogue. For
the case of just two models, Type-I and Type-II errors provide exactly
this analogue --- note that if both errors go to 0 as $n \rightarrow
\infty$, this implies consistency. Thus, the {\em practical\/}
importance of the present work, for us, is mostly that model
comparison by switching defines, like Bayes, a robust null hypothesis
test --- providing Type-I errors irrespective of the stopping rule and
thus more in line with actual practice --- yet has better Type-II
error behaviour, allowing the Type-II error to become small (i.e. the
power to go to $1$) whenever the true distribution sits at a distance
of order $\sqrt{(\log \log n)/n}$ rather than $\sqrt{(\log n)/n}$, as
with Bayes. We only showed robustness  for singleton $\cM_0$,
however, and our simulations show that it may fail for composite $\cM_0$, so {\em the\/} major goal for future work is therefore, to
come up with methods that are  robust to optional stopping
also under composite $\cM_0$.

\paragraph{How special is the switch distribution? }
Since Yang proved that in general, the conflict between consistency
and risk-optimality is not resolvable, one might argue that any model
selection rule just picks some position in the spectrum of behaviours
of consistency vs. risk-optimality. For example, one might have a modified HQ criterion which picks ${\cal M}_1$ if, using the same setup and notation
as in (\ref{eq:AICLIL}), 
\begin{equation}\label{eq:logparty}
\left|  \sum_{i=1}^n \tilde X_i \right| \geq \sqrt{ n \log \log\log n}.
\end{equation}
By the central limit theorem, such a method will be consistent, yet
when combined with an efficient estimator will achieve the minimax
estimation rate up to a $\log \log \log n$ factor, improving on the
switch criterion by an additional logarithm. Note however that both
the switch distribution and HQ (with $c > 1$) achieve {\em strong\/}
consistency. The meaning of strong consistency is illustrated in
Figure~\ref{fig:optstopna} above: it means that, from some $n$ onward,
the wrong model will never be selected any more, no matter how long
one keeps sampling. It is easy to see from the law of the iterated
logarithm that any strongly consistent method can have rate no faster
than order $(\log \log n)/n$ --- in particular, (\ref{eq:logparty}) is
not strongly consistent. Thus, in this sense both switching and HQ do
take a special place in the consistency vs. risk-optimality spectrum
as obtaining the fastest rates compatible with strong consistency,
which may be viewed as asymptotic robustness to optional stopping.
While this may mostly be of theoretical interest, the switch
distribution also takes a special place in terms of its nonasymptotic
robustness to optional stopping: again, the law of the iterated
logarithm implies that any model comparison method that defines a
robust hypothesis test cannot achieve estimation rate better than
order $(\log \log n)/n$. Again, the main open question here is whether one can
modify it so that robustness for composite ${\cal M}_0$ is achieved as well.

\paragraph{Future Work --- Limitations of the Switch Distribution and Our Results}
Whereas the results in this paper all apply to the original switch
distribution as defined by Van Erven et al. (2007) and a
simplification thereof, for full robustness to optional stopping with
composite ${\cal M}_0$, some substantial changes have to be made, as
suggested by the results in Figure~\ref{fig:optstopvoor}.  Initial
research suggests that such a modification of the switch distribution
might indeed be constructed, based on techniques in \cite{RamdasB15};
whereas, compared to Bayes factor testing, in the current switch
criterion, $p_{B,1}$ is modified to another distribution and $p_{B,0}$
can remain the same, in this new version we would also have to change
$p_{B,0}$ --- the resulting distribution would not have a Bayesian
interpretation any more. While this work is still under development, to
avoid the nonrobustness seen in Figure~\ref{fig:optstopvoor} as much
as possible, for the time being we recommend using flat priors (but in
this case, not completely flat - Jeffreys' prior on $\mu$ is improper,
in which case Theorem \ref{thm:optstopnull} holds in none of the
scenarios and simulations --- not reported here --- show that optional
stopping robustness is violated).

Another limitation lies not in the switch distribution, but in our
results: these are restricted to two nested exponential family
models. It would be interesting to extend them to more than two models
--- highlighting the distinction between model selection and testing
--- and going beyond exponential families. We are hopeful that
switching still behaves well in such contexts --- we note that the
risk rate convergence results of Van
Erven et al. (\citeyear{Erven2012}) were for countable,
possibly infinite collections of completely general models --- but
they invariably dealt with the cumulative risk. While all our
experiments suggest that small cumulative risk usually goes together
with small instantaneous risk, formal analysis of the switch
criterion's instantaneous risk is far more difficult, and the present
paper heavily relies on sufficiency to do so --- so extension of our
results beyond exponential families would be difficult.

Before doing so, we would prefer to modify the switch distribution
further, since the present version has a drawback when used in
nonsequential settings: the precise results it gives are dependent on
the order of the data, even if all the models under consideration are
i.i.d. Thus, it would be interesting and challenging to design an
alternative, order-independent method that, like the switch
distribution, is strongly consistent, near rate- and power-optimal,
and is robust to optional stopping under composite $\cM_0$. Such a
method would essentially truly achieve the best of the three worlds we
considered in this paper --- and this is the method we aim for in our
future research.

\section* {Acknowledgements}

The central result of this paper, Theorem~\ref{thm:risk}, already
appeared in the Master's Thesis (Van der Pas,
\citeyear{Vanderpas2013}) for the (very) special case where $m_1 = 1$
and $m_0 = 0$, but the proof supplied there contained a (serious but
repairable) error. We are grateful to Tim van Erven for pointing this
out to us. We are also thankful to the anonymous referees and to
Hannes Leeb for raising the issue of whether the switch distribution
has a `special' place on the spectrum of a model selection criterion's
possible risk and consistency behaviors.

\appendix

We start by listing some well-known properties of
exponential families which we will repeatedly use in the proofs. Then,
in Section~\ref{sec:mainproof}, we provide a sequence of technical
lemmata that lead up to the proof of our main result,
Theorem~\ref{thm:risk}. Finally, in
Section~\ref{app:realswitch}, we compare the switch distribution and criterion as defined
here to the original switch distribution and criterion of Van Erven et al. (\citeyear{Erven2012}).

\paragraph{Additional Notation}
Our results will often involve displays involving several
constants. The following abbreviation proves useful: when we write `for positive
constants $\vec{c}$, we have ...', we mean that there exist some
$(c_1, \ldots, c_N) \in \reals^N$, with $c_1, \ldots, c_N > 0$, such
that ... holds; here $N$ is left unspecified but it will always be clear
from the application what $N$ is.  
Further, for positive constants $\vec{b} = (b_1,b_2,b_3)$, we define $\restje$ as
$$
\restje = \begin{cases} 1 & \text{if $n < b_1$} \\
b_2 e^{-b_3 n} & \text{if $n \geq b_1$},
\end{cases}
$$
and we frequently use the following fact. Suppose that $\cE_1, \cE_2, \ldots$ is a
  sequence of events such that $\prob(\cE_n)
\leq \restje$. Then we also have, for any event $\cA$, and for all $n$,
\begin{equation}\label{eq:vhp}
  \prob(\cA, \cE^c_n) \geq 
  \prob(\cA) - \restje,
\end{equation}
as is immediate from $\prob(\cA, \cE^c_n) =
\prob(\cA) - \prob(\cA, {\cE}_n) \geq
\prob(\cA) - \prob({\cE}_n)$.

The components of a vector $\mu \in \mathbb{R}^n$ are given by $(\mu_1, \mu_2, \ldots, \mu_{n})$. If the vector already has an index, we add a comma, for example $\mu_1 = (\mu_{1,1}, \mu_{1,2}, \ldots, \mu_{1,n})$. A sequence of vectors is denoted by $\mu^{(1)}, \mu^{(2)}, \ldots$.

\section{Definitions Concerning and Properties of Exponential Families}
\label{app:expfam}

The following definitions and properties can all be found in the standard reference
\citep{BarndorffNielsen1978} and, less formally, in \citep[Chapters
18 and 19]{Grunwald2007}.

A $k$-dimensional exponential family is a set of distributions on
$\mathcal{X}$, which we invariably represent by the corresponding set
of densities $\{p_\theta \mid \theta \in \Theta\}$, where $\Theta \subset
\mathbb{R}^k$, such that any member $p_\theta$ can be written as
\begin{equation}\label{eq:defexp}
p_\theta(x) = \frac{1}{z(\theta)}e^{\theta^T\phi(x)}r(x) = e^{\theta^T\phi(x) - \psi(\theta)}r(x) ,
\end{equation}
where $\phi(x) = (\phi_1(x), \ldots, \phi_k(x))$ is a \emph{sufficient statistic}, $r$ is a non-negative function called the \emph{carrier}, $z$
the \emph{partition function} and $\psi(\theta) = \log{z(\theta)}$. We assume the representation \eqref{app:expfam} to be \emph{minimal}, meaning that the components of $\phi(x)$ are linearly independent.

The parameterization in (\ref{eq:defexp}) is referred to as the
\emph{canonical} or \emph{natural parameterization\/}; we only
consider families for which the set $\Theta$ is open and connected.
Every exponential family can alternatively be parameterized in terms
of its \emph{mean-value parameterization}, where the family is
parameterized by the mean $\mu = \E_{\theta}[\phi(X)]$, with $\mu$
taking values in $M \subset \mathbb{R}$, where $\mu$ as a function of
$\theta$ is smooth and strictly increasing; as a consequence, the set
$M$ of mean-value parameters corresponding to an open and connected
set $\Theta$ is itself also    open and connected. Whenever for data
$x_1, \ldots, x_n$, we have $\frac{1}{n}\sum_{i=1}^n \phi(x_i) \in M$,
then the maximum likelihood is uniquely achieved by the $\mu$ that is
itself equal to this value,
\begin{equation}\label{eq:mle}
\widehat\mu(x^n) =
\frac{1}{n}\sum_{i=1}^n \phi(x_i).
\end{equation}
We thus define the maximum likelihood estimator (MLE) to be equal to
(\ref{eq:mle}) whenever
\begin{equation}\label{eq:meanisml}
\frac{1}{n} \sum_{i=1}^n  \phi(X_i) \in
M. 
\end{equation} 
Since the result below which directly involves the MLE
(Lemma~\ref{lem:8.1}) does not depend on its value for $x^n$ with
$\frac{1}{n}\sum_{i=1}^n \phi(x_i) \not \in M$, we can leave
$\widehat{\mu}(x^n)$ undefined for such values.  However, if we want
to use the MLE as a `sufficiently efficient' estimator as used in the
statement of Theorem~\ref{thm:risk}, we need to define
$\widehat{\mu}(x^n)$ for such values in such a way that the
`sufficiently efficient property' (\ref{eq:suffeff}) is satisfied. The
following examples show various ways of constructing such sufficiently
efficient estimators.

\begin{example}{\bf [Sufficient Efficiency for MLE's for squared
    (standardized) error and Hellinger]\ }{\label{ex:mle}\rm 
For many full families such as the full
(multivariate) Gaussians, Gamma and many others, (\ref{eq:meanisml})
holds $\mu$-almost surely for each $n$, for all $\mu \in M$. 
If we compare two families ${\cal M}_0$ and ${\cal M}_1$  given in their mean-value
parameterization with $M_0 \subset M_1$  where ${\cal M}_1$ is any such family,
then
the MLE is almost surely well-defined for $M_1$ and thus we need not worry about the issue indicated above. We can then  take
$\breve{\mu}_1:= \widehat\mu_1$ to be the MLE for $\cM_1$. To get a sufficiently efficient estimator for $M_0$, we take $\breve{\mu}_0$ to be the  projection of $\widehat\mu_1$ 
on the first $m_0$ coordinates
(usually (\ref{eq:meanisml}) will still hold for $\cM_0$ and then
this $\breve\mu_0$ will also be the MLE for $\cM_0$). This pair of estimators
will be sufficiently efficient for (standardized) squared error and
squared Hellinger distance, i.e. (\ref{eq:suffeff}) holds for these three losses. To show
this, note that from Proposition~\ref{prop:withinconstants}, Eq. (\ref{eq:seconddisp}), we see that
it is sufficient to show that (\ref{eq:suffeff}) holds for the squared
error loss. Since the $j$-th component of $\widehat\mu_1$ is equal to
$n^{-1} \sum_{i=1}^n \phi_j(X_i)$ and $\E_{\mu_1}[n^{-1} \sum_{i=1}^n
\phi_j(X_i) ] = \mu_{1,j}$ and 
$
\text{\sc var}_{\mu_1} \left[n^{-1}
\sum_{i=1}^n \phi_j(X_i)\right] = n^{-1} \text{\sc var}_{\mu_1} \left[
\phi_j(X_1)\right],
$
it suffices to show that
$$
\sup_{\mu_1 \in M'_1} \sup_{j = 1, \ldots, m_1} \text{\sc var}_{\mu_1} \left[
\phi_j(X_1)\right] = O\left(1  \right),
$$
which is indeed the case since $M'_1$ is a CINECSI set, so that the
variance of all $\phi_j$'s is uniformly bounded on $M'_1$
\citep{BarndorffNielsen1978}.}
\end{example}
\begin{example}{\bf [Other sufficiently efficient estimators for
    squared (standardized) error and Hellinger]\ }{\rm For models such
    as the Bernoulli or multinomial, (\ref{eq:meanisml}) may fail to
    hold with positive probability: the full Bernoulli exponential
    family does not contain the distributions with $P(X_1=1) = 1$ and
    $P(X_1= 0)=1$, so if after $n$ examples, only zeros or only ones
    have been observed, the MLE is undefined. We can then go either of
    three ways. The first way, which we shall not pursue in detail
    here, is to work with so-called `aggregate' exponential families,
    which are extensions of full families to their limit points. For
    models with finite support (such as the multinomial) these are
    well-defined \cite[page 154--158]{BarndorffNielsen1978} and then
    the MLE's for these extended families are almost surely
    well-defined again, and the MLE's are sufficiently efficient by
    the same reasoning as above. Another approach that works in some
    cases (e.g. multinomial) is to take $\breve{\mu}_1$ to be a
    truncated MLE, that, at sample size $n$, maps $X^n$ to the MLE
    within some CINECSI subset $M^{(n)}_1$ of $M_1$, where $M_1^{(n)}$
    converges to $M_1$ as $n$ increases in the sense that $\sup_{\mu
      \in M_1^{(n)},\mu' \in M_1 \setminus M_1^{(n)}} \| \mu - \mu'
    \|^2_2 = O(1/n)$. The resulting truncated MLE, and its projection
    on $M_0$ (usually itself a truncated MLE) will then again be
    sufficiently efficient. This approach also works if the models
    $\cM_0$ and $\cM_1$ are not full but restricted families to begin
    with. For full families though, a more elegant approach than
    truncating MLE's is to work with Bayesian posterior MAP estimates
    with conjugate priors. For steep exponential families (nearly all
    families one encounters in practice are steep), one can always
    find conjugate priors such that the Bayes MAP estimates based on
    these priors exist and take a value in $M_1$ almost surely
    \citep{GrunwaldD05}. They then take the form $\breve{\mu}_1 =
    \sum_{i=1}^n (\phi(X_i) + \lambda_0 \mu^\circ_1) / (n+
    \lambda_0)$, where $\lambda_0 > 0$ and $\mu^\circ_1 \in M_1$ are
    determined by the prior. $\breve\mu_0$ can then again be taken to
    be the projection of $\breve\mu_1$ onto $M_0$. Under the
    assumption that $\mu_1$ is contained in a CINECSI set $M'_1$, one
    can now again show, using the same arguments as in
    Example~\ref{ex:mle}, that such estimators are sufficiently
    efficient for squared (standardized) error and Hellinger loss.}
\end{example}
\begin{example}{\bf [Sufficient Efficiency for R\'enyi and KL
    divergence]\ }{\rm As is well-known, for the multivariate Gaussian
    model with fixed covariance matrix, the squared error risk and KL
    divergence are identical up to constant factors, so the
    unrestricted MLE's will still be sufficiently efficient for KL
    divergence. For other models, though, the MLE will not always be
    sufficiently efficient. For example, with the Bernoulli model and
    other models with finite support, to make the unrestricted MLE's
    well-defined, we would have to extend the family to its boundary
    points as indicated in Example~\ref{ex:mle}.  Since, however, for
    any $0 < \mu < 1$ and $\mu'=0$, the KL divergence $D(\mu \| \mu')
    = \infty$ and $\prob_{\mu}(\hat\mu(X^n) = \mu') > 0$, the
    unrestricted MLE in the full Bernoulli model including the
    boundaries will have infinite risk and thus will not be sufficiently
    efficient. The MAP estimators tend to behave better though:
    \cite{GrunwaldD05} implicitly show that for 1-dimensional
    families, under weak conditions on the family (Condition 1
    underneath Theorem 1 in their paper) --- which were shown to hold
    for a number of families such as Bernoulli, Poisson, geometric ---
    sufficient efficiency for the KL divergence still holds for MAP
    estimators of the form above. We conjecture that a similar result
    can be shown for multidimensional families, but will not attempt
    to do so here.}
\end{example}

A standard property of exponential families says that, for any $\mu \in
M$, any distribution ${\mathbb Q}$ on ${\cal X}$ with $\E_{X \sim
  {\mathbb Q}} [\phi(X)] = \mu$, any $\mu' \in M$, we have
\begin{equation}
\E_{X \sim
  {\mathbb Q}} \left[\log \frac{p_{\mu}(X)}{p_{\mu'}(X)} \right] =
\E_{X \sim
  {\mathbb P_{\mu}}} \left[\log \frac{p_{\mu}(X)}{p_{\mu'}(X)} \right] =
D(\mu \| \mu'), \label{eq:robustnesspre}
\end{equation}
the final equality being just the definition of $D(\cdot \|
\cdot)$. Now fix an arbitry sample $x^n$. 
By taking ${\mathbb Q}$ to be the empirical distribution on
${\cal X}$ corresponding to sample $x^n$, it follows from
(\ref{eq:robustnesspre}) that if $\widehat\mu(x^n) \in M$ then also the
following relationship holds for any $\mu'\in M$: 
\begin{equation}
\frac{1}{n} \log{\frac{p_{\widehat \mu(x^n)}(x^n)}{p_{\mu'}(x^n)} } 
= D(\widehat\mu(x^n) \| \mu') \label{eq:robustness}.
\end{equation}
(\ref{eq:robustnesspre}) and (\ref{eq:robustness}) are a direct consequence of the sufficiency of
$\widehat{\mu}_1(X^n)$, and folklore among information theorists. For
a proof of (\ref{eq:robustnesspre}) and more details on
(\ref{eq:robustness}), see e.g. \cite[Chapter 19]{Grunwald2007}, who
calls this the \emph{robustness property} of the KL divergence for
exponential families.

We are now in a position to prove Proposition~\ref{prop:withinconstants}, which we repeat for convenience. 
\paragraph{Proposition~\ref{prop:withinconstants}}
  Let $M$, a product of open intervals, be the mean-value parameter
  space of an exponential family, and let $M'$ be a CINECSI
  subset of $M$. Then there exist positive constants $\vec{c}$  such that for all $\mu, \mu' \in M'$,
\begin{equation}\label{eq:firstdisp}
c_1 \| \mu' - \mu  \|^2_2 \leq c_2 \cdot \std(\mu' \| \mu)
\leq \hellinger(\mu',\mu) \leq \renyi(\mu', \mu)
\leq D(\mu' \| \mu) \leq c_3 \| \mu' - \mu  \|^2_2.
\end{equation}
and for all $\mu' \in M', \mu \in M$ (i.e. $\mu$ is now not restricted to lie in $M'$),
\begin{equation}\label{eq:seconddisp}
\hellinger(\mu',\mu) \leq c_4 \| \mu' - \mu  \|^2_2 \leq c_5 \cdot \std(\mu' \| \mu)
\leq c_6 \| \mu' - \mu  \|^2_2.
\end{equation}
\begin{proof} 
We start with (\ref{eq:firstdisp}).
The third and fourth inequality are immediate by using
  $-\log x \geq 1-x$ 
and Jensen's inequality, respectively. From standard properties of Fisher
  information for exponential families \citep{BarndorffNielsen1978} we
  have that, for any CINECSI (hence compact and bounded away from the
  boundaries of $M$) subset $M'$ of $M$, there exists positive $\vec{C}$ with
\begin{equation}\label{eq:boundedfisher}
0 <  {C}_1 =  \inf_{\mu \in M'} \det I(\mu) < 
\sup_{\mu \in M'} \det I(\mu) = C_2 < \infty,
\end{equation}
from which we infer that for all $\mu' \in M'$, $\mu,\mu'' \in \reals^m$,
\begin{equation}\label{eq:eigen}
C_3 \| \mu- \mu'' \|^2_2 \leq (\mu - \mu'')^T I(\mu') (\mu - \mu'') \leq  C_4 \| \mu- \mu'' \|^2_2,
\end{equation}
for some $0 < C_3 \leq C_4 < \infty$. 
Using (\ref{eq:eigen}), the first inequality is immediate, and the final
inequality  follows straightforwardly from a second-order Taylor
approximation of KL divergence as in
\cite[Chapter 4]{Grunwald2007}. It only remains to establish the second inequality. Now, since $M'$ is CINECSI and
hence compact the fifth (rightmost)  inequality implies 
that there is a $C_5< \infty$ such that $\sup_{\mu,\mu' \in M'} D(\mu' \|
\mu) < C_5$ and hence, via the fourth inequality, that $\sup_{\mu,\mu' \in M'} \renyi(\mu', \mu) < C_5$.  Equality (\ref{eq:rh}) now implies that there is a $C_6$ such
that 
\begin{equation}\label{eq:rianne} \sup_{\mu,\mu' \in M'} \renyi(\mu', \mu) /  \hellinger(\mu',\mu)  < C_6.\end{equation}
Using again (\ref{eq:boundedfisher}), a
second order Taylor approximation as in Van Erven and Harremo\"es (\citeyear{Erven2014}) now gives that
for some constant $C_7 > 0$, $\| \mu - \mu'  \|^2_2 \leq C_7 \renyi(\mu', \mu)$
for all $\mu,\mu' \in M'$. The first result, (\ref{eq:firstdisp}), now follows upon combining this
with (\ref{eq:rianne}).  

As to (\ref{eq:seconddisp}), the second and third inequality are
immediate from (\ref{eq:eigen}). For the first inequality, note that,
since $M'$ is CINECSI and we assume $M$ to be a product of open
intervals, there must exist another CINECSI subset $M''$ of $M$
strictly containing $M'$ such that $\inf_{\mu' \in M', \mu \in M
  \setminus M''} \| \mu' - \mu \|_2^2 = \delta$ for some $\delta > 0$.
We now distinguish between $\mu$ in (\ref{eq:seconddisp}) being an
element of (a) $M''$ or (b) $M \setminus M''$. For case (a)
(\ref{eq:firstdisp}), with $M''$ in the role of $M'$, gives that there
is a constant $C_8$ such that for all $\mu \in M''$,
$\hellinger(\mu',\mu) \leq C_8 \| \mu'- \mu \|_2^2$. For case (b), $\mu
\in M \setminus M''$, we have $\| \mu'- \mu \|^2_2 \geq \delta$ and,
using that squared Hellinger distance for any pair of distributions is
bounded by $2$, we have $\hellinger(\mu',\mu) \leq (2/ \delta) \|
\mu'- \mu \|_2^2$. Thus, by taking $c_4 = \max \{C_8, 2/\delta\}$, case
(a) and (b) together establish the first inequality in
(\ref{eq:seconddisp}).
\end{proof}

\section{Preparation for Proof of Main Result: Results on Large
  Deviations}
 Let $\cM_1$ and $M_1$ be as in
Theorem~\ref{thm:risk}. For the following result,
Lemma~\ref{lem:klbound}, we set $\widehat\mu'_1(X^n) := n^{-1} \sum
\phi(X_i)$, so that $\widehat\mu'_1(X^n) = \widehat\mu_1(X^n)$
whenever $n^{-1} \sum \phi(X_i) \in M_1$. It is essentially a
multidimensional extension of a standard information-theoretic result,
with KL divergence replaced by squared error loss. This standard
result states the following: whenever ${\cal M}_1$ is a
single-parameter exponential family (that is, $m_1 =1$), then for any
$\mu \in M_1$, all $a,a'> 0$ with $\mu+ a \in M_1$, $\mu - a' \in
M_1$,
\begin{equation}\label{eq:klineq}
\prob_{\mu}(\widehat\mu'_1(X^n) \geq \mu + a) \leq e^{-n D(\mu +a \| \mu)}.
\ \ ;\ \ \prob_{\mu}(\widehat\mu'_1(X^n) \leq \mu - a' ) \leq e^{-n D(\mu -a' \| \mu)}.
\end{equation}
For a simple proof, see \citep[Section 19.4.2]{Grunwald2007}; for
discussion see \citep{Csiszar1984} --- the latter reference gives a
multidimensional extension of (\ref{eq:klineq}) but of a very
different kind than Lemma~\ref{lem:klbound} below. To prepare for the
lemma, let $\cM_1$ and $M_1$ be as in Theorem~\ref{thm:risk} and, for
any $\mu \in M_1$ and any $\vec{a},\vec{b} \in
\mathbb{R}^{m_1}_{>0}$, define the $\ell_{\infty}$-rectangle
$R_{\infty}(\mu, \vec{a},\vec{b}) = \{ \mu' \in \reals^{m_1} :
\forall j = 1, \ldots, m_1, - b_j \leq \mu_j' - \mu_j \leq a_j \}$.
\begin{lemma}\label{lem:klbound}
Let $\cM_1$ and $M_1$ be as in Theorem~\ref{thm:risk} and fix an
arbitrary CINECSI subset $M'_1$ of $M_1$. Then there is a
$c > 0$ (depending on $M'_1$) such that, for all $\mu \in M_1$, all $n$, all $\vec{a},\vec{b} \in
\mathbb{R}^{m_1}_{>0}$ such that $R_{\infty}(\mu, \vec{a},\vec{b}) \subset M'_1$,
\begin{equation}\label{eq:klineqb}
\prob_{\mu}( \widehat\mu'_1(X^n) \not \in R_{\infty}(\mu,\vec{a},\vec{b}))  \leq 2 m_1 
e^{-n  c\cdot (\min_j \min \{a_j, b_j\})^2}.
\end{equation}
\end{lemma}
\begin{proof}
For $j =1, \ldots, m_1$, $d \in \reals$, let $\vec{e_j}$ represent the $j$th standard basis vector,
  such that $\mu + d\vec{e_j} = (\mu_1, \ldots, \mu_{j-1}, \mu_j+d,
  \mu_{j+1}, \ldots, \mu_{m_1})$, and 
  let $D_{\mu + d\vec{e_j}} := D(\mu + d\vec{e_j} \| \mu)$.  We now have that there exist
  constants  $c_{a,1}, \ldots, c_{a,m_1}, c_{b,1}, \ldots, c_{b,m_1} > 0$ such
  that for \break $c:=  \min \{c_{a,1}, \ldots, c_{a, m_1},  c_{b,1}, \ldots, c_{b,m_1} \}$, all $n$, 
\begin{align*}
\prob_{\mu}(\widehat\mu_1(X^n) \not \in R_{\infty}(\mu,\vec{a}, \vec{b}))  & \leq
\sum_{j=1}^{m_1} \prob_{\mu}(\widehat\mu_{1,j}(X_n) \geq \mu_j + a_j) + \sum_{j=1}^{m_1} \prob_{\mu}(\widehat\mu_{1,j}(X^n) \leq \mu_j - b_j) \\
& \leq \sum_{j=1}^{m_1} \left( e^{- n D_{\mu + a_j\vec{e_j}}} + e^{- n D_{\mu - b_j\vec{e_j}}} \right) 
\leq \sum_{j=1}^{m_1} \left( e^{- n c_{a,j} a_j^2} + e^{- n c_{b,j} b_j^2}
\right)  \\ & \leq  2 m_1 
e^{-n  c\cdot (\min_j \min \{a_j,b_j\})^2 },
\end{align*}
Here the first inequality follows from
the union bound, and the second follows by applying, for each of the
$2 m_1$ terms, (\ref{eq:klineq}) above to the one-dimensional
exponential sub-family
$\{ p_{\mu} \mid \mu \in M_1 \cap \{ \mu: \mu = \mu + d\vec{e_j} \text{\ for
  some $d \in {\mathbb R}$} \}\}$.  The third follows by
Proposition~\ref{prop:withinconstants} together with the equivalence
of the $\ell_2$ and sup norms on $\mathbb{R}^{m_1}$, and the final
inequality is immediate. 
\end{proof}

\begin{lemma} \label{lem:likbound} Under conditions and notations as in
  Theorem~\ref{thm:risk}, let $\mu,\mu'$ be  elements of $M_1$ and suppose
  $X^N = (X_{n_1}, \ldots, X_{n_2})$ is a sequence of i.i.d.
  observations of length $N$ from $p_{\mu}$. 
  Then, for any $A \in \mathbb{R}$:
\begin{equation}\label{eq:renyibound}
\prob_{\mu}\left(\log{ \frac{p_{\mu}(X^N)
    }{p_{\mu'}(X^N)}} < A  \right) \leq
e^{\frac{1}{2}A}e^{-\frac{N}{2} \renyi(\mu', \mu)}.
\end{equation}  
\end{lemma}
\begin{proof}
For any $A$, by Markov's inequality:
\begin{align}
\prob_{\mu}&\left(\log{ \frac{p_{\mu}(X^N) }{p_{\mu'}(X^N)}} < A  \right) = \prob_{\mu}\left(\left( \frac{p_{\mu'}(X^N) }{p_{\mu}(X^N )} \right)^{\frac{1}{2}} > e^{-\frac{1}{2}A}  \right) 
\leq e^{\frac{1}{2}A}\E_{\mu}\left[ \left( \frac{p_{\mu'}(X^N) }{p_{\mu}(X^N )} \right)^{\frac{1}{2}}\right] \notag \\
&= e^{\frac{1}{2}A}\left(\E_{\mu}\left[ \left( \frac{p_{\mu'}(X_{n_1}) }{p_{\mu}(X_{n_1})} \right)^{\frac{1}{2}}\right] \right)^{N}
= e^{\frac{1}{2}A}e^{\log{\left(\E_{\mu}\left[ \left( \frac{p_{\mu'}(X_{n_1}) }{p_{\mu}(X_{n_1})} \right)^{\frac{1}{2}}\right] \right)^{N} } }\notag\\
&= e^{\frac{1}{2}A}e^{-\frac{N}{2}\left( - \frac{1}{1-1/2}\log{\E_\mu\left[ \left( \frac{p_{\mu'}(X_{n_1}) }{p_{\mu}(X_{n_1})} \right)^{\frac{1}{2}}\right]  } \right) }
=  e^{\frac{1}{2}A}e^{-\frac{N}{2} \renyi(\mu, \mu')} \label{eq:renyi}.
\end{align}
\end{proof}

\begin{proposition}\label{prop:niceml}
Let $\cM_0, \cM_1, M_0, M_1$ be as in Theorem~\ref{thm:risk} and let
$M'_1$ be a  CINECSI subset of $M_1$. Then there exists another, larger, 
CINECSI subset $M''_1$ of $M_1$ and positive constants $\vec{b}$ 
such that $M'_1$ is itself a CINECSI
subset of $M''_1$ and  for both $j \in \{0,1\}$, the ML estimator $\widehat\mu_j(x^n)$ satisfies 
$$
\sup_{\mu \in M'_1} P_{\mu}(\widehat{\mu}_j(X^n) \not \in M''_1) \leq \restje.
$$
\end{proposition}
\begin{proof}
$M_1$ can be written as in (\ref{eq:parameters}), and hence we can
define a set 
$$M''_1 =
[\zeta_{1,1}^{*}, \eta_{1,1}^{*}] \times \ldots
\times [\zeta_{1, m_1}^{*}, \eta_{1, m_1}^{*}]$$ for
values $\zeta_{1,j}^{*}, \eta_{1,j}^{*} \in \reals$
such that $M''_1$ is a CINECSI subset of $M_1$. Since $M'_1$ is
connected with compact closure in interior of $M_1$ and $M''_1$ is a
subset of $M_1$, we can choose the $\zeta
_{1,j}^{*}, \eta_{1,j}^{*} \in \reals$ such that $M'_1$ is itself
a CINECSI subset of $M''_1$. Since $M'_1$ is connected and its
closure is in the interior of $M''_1$ which is itself compact, it
follows that there is some $\delta> 0$ such that, for all $\mu'_1 \in
M'_1, \mu''_1 \not \in M''_1$, all $j \in \{1, \ldots, m_1\}$, it holds
$|\mu'_{1,j} - \mu''_{1,j}| > \delta$. It now follows from
Lemma~\ref{lem:klbound}, applied with $\vec{a}$ chosen such that
$R_{\infty}(\mu',\vec{a}) = M''_1$,   that for every $\mu'\in M'_1$,
all $n$,
$$
\prob_{\mu'}\left(
\widehat{\mu}_{1}(X^n) \not \in M_1'' \right) \leq C_1 e^{- n C_2 \delta^2}
$$
for some constants $C_1, C_2$. Here we used that by construction, each entry of $\vec{a}$
must be at least as large as $\delta$. Since $\widehat{\mu}_{1,j}(x^n)$ and
$\widehat{\mu}_{0,j}(x^n)$ coincide for $0 < j \leq m_0$  and
$\widehat{\mu}_{0,j}(x^n)$ is constant for $m_0 < j \leq m_1$, the result
follows for $\widehat{\mu}_0(x^n)$ as well.
\end{proof}

\section{Preparation for Proof of Main Result: Results on Bayes
  Factor Model Selection}
 \begin{lemma} \label{lem:8.1} Let $\cM_0, \cM_1, M_0, M_1$ be as in
   Theorem~\ref{thm:risk} and let, for $j \in \{0,1\}$, $M'_j$ be a
   CINECSI subset of $M_j$. For both $j \in \{0,1\}$, there exist
   positive constants $\vec{c}, \vec{b}$ such that for all $\mu_1 \in M'_1$,
\begin{equation}\label{eq:laplace}
c_1 \leq {n}^{-m_j/2} \cdot \frac{ p_{\widehat\mu_j(X^n)}(X^n)}{p_{B,j}(X^n)} \leq c_2,
\end{equation}
with $\prob_{\mu_1}$-probability at least $1- \restje$. 
\end{lemma}
\begin{proof} 
For a Bayesian marginal distribution $p_B$ defined relative to
$m$-dimensional exponential family $\cM$ given in its mean-value
parameterization $M$, with a prior $\omega(\cdot)$ that is continuous and strictly positive
  on $M$, we have as a consequence of the familiar Laplace
  approximation of the Bayesian marginal distribution of exponential
  famlies as in e.g. \citep{Kass1995},
\begin{equation*}
p_B(x^n) \sim \left(\frac{n}{{2\pi}}\right)^{- m/2} \cdot \frac{\omega(\widehat\mu(x^n))}{\sqrt{\det I(\widehat\mu(x^n))}} p_{\widehat\mu(x^n)}(x^n).
\end{equation*}
As shown in Theorem 8.1 in \citep{Grunwald2007}, this statement holds
uniformly for all sequences $x^n$ with ML estimators in any fixed {\em
  CINECSI\/} subset $M'$ of $M$.  By compactness of $M'$, and by
positive definiteness and continuity of Fisher information for
exponential families, the quantity $\omega(\widehat\mu) / \sqrt{ \det
  I(\widehat\mu)}$ will be bounded away from zero and infinity on such
sequences, and, applying the result to both the families $\cM_0$ and
$\cM_1$ it follows that there exist $c_1, c_2 > 0$ such that for all
$n$ larger than some $n_0$, uniformly for all sequences $x^n$ with
$\widehat\mu_j(x^n) \in M'_j$, we have:
\begin{equation}\label{eq:laplaceb}
c_1 \leq {n}^{-m_j/2} \cdot \frac{
  p_{\widehat\mu_j(x^n)}(x^n)}{p_{B,j}(x^n)} \leq c_2. 
\end{equation}
The result now follows by combining this statement with Proposition~\ref{prop:niceml}.
\end{proof}

\begin{lemma}\label{lem:bayesfactorrenyi}
Let $\cM_0, \cM_1, M_0$, $M_1$ and the Bayesian marginal
distribution $p_{B,0}$ be as in Theorem~\ref{thm:risk}.
Let $M'_1$ be a CINECSI subset of $M_1$. 
Then there exist positive constants $\vec{c}$ and $\vec{b}$ such
that for all $n$, all $\mu_1 \in M'_1$, all $A \in {\mathbb R}$,
$$
\prob_{\mu_1}\left(\log{ \frac{p_{B,1}(X^n)
    }{p_{B,0}( X^n )}} < A  \right) \leq
n^{m_1/2} \cdot c_1 \cdot e^{\frac{1}{2} c_2 A}e^{-\frac{n}{2} c_3 \|
  \mu_1 - \mu_0 \|_2^2} + \restje, 
$$
where for each $\mu_1$, $\mu_0 = \proj(\mu_1)$ as in
(\ref{eq:mu0b}).
\end{lemma}

\begin{proof} 
  Fix constants $C_1, C_2$ such that they are smaller and larger respectively than the constants $c_1, c_2$ from Lemma~\ref{lem:8.1} and define $$\cE_n= \left\{ X^n: C_1 \leq
    n^{-m_1/2} \frac{p_{\widehat\mu_1(X^n)}(X^n)}{p_{B,1}(X^n)} \leq C_2 \right\}.$$ Using
    Lemma~\ref{lem:8.1}, we have that there exists positive $\vec{b}$
    such that for all $A \in \reals$, 
\begin{align}\label{eq:start2}
& \prob_{\mu_1}\left(\log{ \frac{p_{B,1}(X^n)
    }{p_{B,0}( X^n )}} < A  \right) \notag \\
= & \prob_{\mu_1}\left(\log{ \frac{p_{B,1}(X^n)
    }{p_{B,0}( X^n )}} < A , \cE_n \right) + \prob_{\mu_1}\left(\log{ \frac{p_{B,1}(X^n)
    }{p_{B,0}( X^n )}} < A , {\cE}^c_n \right)  \notag \\
\leq & \prob_{\mu_1}\left(\log{ \frac{ C_2^{-1}n^{-m_1/2}  p_{\widehat{\mu}_1(X^n)}(X^n)
    }{p_{B,0}( X^n )}} < A , \cE_n \right)
+ 
\restje \notag \\
\leq  & \prob_{\mu_1}\left(\log{ \frac{C_2^{-1}n^{-m_1/2}   p_{{\mu}_1}(X^n)
    }{p_{B,0}( X^n )}} < A  \right)
+ 
\restje \notag \\ & =  \prob_{\mu_1}\left(\log{ \frac{ p_{{\mu}_1}(X^n)
    }{p_{B,0}( X^n )}} < A + \log C_2 n^{m_1/2}  \right)
+ \restje.
\end{align}
To bound this probability further, we need to relate $p_{B,0}$ to
$p_{B',0}$, the Bayesian marginal likelihood under model $M_0$ under a
prior with support restricted to a compact set $M'_0$. To define
$M'_0$, note first that there must exist a CINECSI subset, say
$M''_1$, of $M_1$ such that $M'_1$ is itself a CINECSI subset of
$M''_1$. Take any such $M''_1$ and let $M'_0$ be the closure of $M''_1
\cap M_0$. Given $\omega$,  the prior density on $\Pi'(M_0)$ used in the
definition of $p_{B,0}$, define $\omega'(\nu) = \omega(\nu)/\int_{\nu
  \in \Pi'(M'_0)} \omega(\nu) d \nu$ as the prior density restricted to
and normalized on $\Pi'(M'_0)$  and let $p_{B',0}$ be the corresponding
Bayesian marginal density on $X^n$. 

To continue bounding (\ref{eq:start2}), define 
$$\cE'_n= \left\{ X^n: C_3 \leq 
 n^{-m_0/2} \frac{p_{\widehat\mu_0(X^n)}(X^n)}{p_{B,0}(X^n)} \leq C_4   \text{\ and \ }
C_3 \leq 
 n^{-m_0/2} \frac{p_{\widehat\mu_0(X^n)}(X^n)}{p_{B',0}(X^n)} \leq C_4  \right\},
$$
with $C_3$ and $C_4$ smaller and larger respectively than the
constants $c_1$ and $c_2$ resulting from Lemma~\ref{lem:8.1} (note that Lemma~\ref{lem:8.1} can
be applied to $p_{B',0}$ as well, by taking  $M_0$ in that
lemma to be the interior of $M'_0$ as defined here). 
Set $C_5 > C_4/C_3$, and note that for any $A_1 \in \reals$, abbreviating
$\prob_{\mu_1}\left(\log{ \frac{p_{\mu_1}(X^n) }{C_5 p_{B',0}(X^n)}} <
  A_1\right)$ to $p^*$, we have
\begin{align}\label{eq:boting} & 
\prob_{\mu_1}\left(\log{ \frac{p_{\mu_1}(X^n)
    }{p_{B,0}(X^n)}} < A_1 \right) \notag \\ = &
\prob_{\mu_1}\left(\log{ \frac{p_{\mu_1}(X^n)
    }{p_{B,0}(X^n)}} < A_1 ,\  \frac{p_{B_0}(X^n)}{p_{B',0}(X^n)} < C_5  \right)
+ \prob_{\mu_1}\left(\log{ \frac{p_{\mu_1}(X^n)
    }{p_{B,0}(X^n)}} <A_1 ,\  \frac{p_{B_0}(X^n)}{p_{B',0}(X^n)} \geq C_5  \right)
\notag \\  \leq & 
\prob_{\mu_1}\left(\log{ \frac{p_{\mu_1}(X^n)
    }{C_5 p_{B',0}(X^n)}} <A_1  \right)
+ \prob_{\mu_1}\left(p_{B,0}(X^n) \geq C_5 p_{B',0}(X^n) \right)
\notag \\ = & p^*
+ \prob_{\mu_1}\left(p_{B,0}(X^n) \geq C_5 p_{B',0}(X^n) \right)
\notag \\ \leq & p^* + \prob_{\mu_1}\left(\frac{p_{B,0}(X^n)}{p_{B',0}(X^n)} \geq C_5  , \cE'_n \right)
+ \prob_{\mu_1}\left(\frac{p_{B,0}(X^n)}{p_{B',0}(X^n)} \geq C_5  ,({\cE}'_n)^c
\right) \notag \\ \leq & p^* + 0 + \restje.
\end{align} 

Now it only remains to bound $p^*$. To this end, let
\begin{equation}\label{eq:priorcond}
C_6 := \int_{\nu \in \Pi'(M'_0)} \sqrt{\omega(\nu)} d \nu.
\end{equation} 
Since $M'_0$ has compact closure in the interior of $M_0$ and we are
assuming that $\omega$ has full support on $M_0$, we have that $C_6 <
\infty$.

Now using Markov's inequality as in the proof of Lemma~\ref{lem:likbound},
that is, the first line of (\ref{eq:renyi})
with $p_{B',0}$ in the role of $p_{\mu'}$, gives, for any $A_2 \in \reals$,
\begin{equation}\label{eq:hoting}
\prob_{\mu_1}\left(
\log \frac{p_{\mu_1}(X^n)
    }{p_{B',0}(X^n)} < A_2  \right) \leq e^{\frac{1}{2}{A_2}}
\E_{\mu_1}\left[ \left( \frac{p_{B',0}(X^n) }{p_{\mu_1}(X^n)} \right)^{\frac{1}{2}}\right].
\end{equation}
The expectation on the right can be further bounded, defining $\omega''
= \sqrt{\omega}/C_6$ and noting that $\omega''$ is a probability
density, as
\begin{align*}
\E_{\mu_1}\left[ \left( \frac{p_{B',0}(X^n) }{p_{\mu_1}(X^n)} \right)^{\frac{1}{2}}\right]
& \leq 
\E_{\mu_1}\left[ \left( \frac{\int_{\nu \in \Pi'(M'_0)} \omega(\nu)^{1/2}
      p_{\nu}(X^n)^{1/2} d \nu}{p_{\mu_1}(X^n)^{1/2}} 
\right)
\right] \\ & = 
C_6 \cdot \E_{\mu \sim \omega''} \E_{\mu_1} \left[ \left(
  \frac{p_{\mu}(X^n) }{p_{\mu_1}(X^n)}\right)^{\frac{1}{2}}
\right] \leq  C_6 \cdot\E_{\mu_1} \left[ \left(
  \frac{p_{\mu^\circ}(X^n) }{p_{\mu_1}(X^n)}\right)^{\frac{1}{2}}
\right],
\end{align*}
where $\mu^\circ \in M'_0$ achieves the supremum of $E_{\mu_1} \left[ \left(
  \frac{p_{\mu^\circ}(X^n) }{p_{\mu_1}(X^n)}\right)^{\frac{1}{2}}
\right]$ within $M'_0$. By compactness of $M'_0$ and continuity, this
supremum is achieved. The final term can be rewritten, following the same steps
as in the second and third line of (\ref{eq:renyi}), as
\begin{equation}\label{eq:joting}
\E_{\mu_1} \left[ \left(
  \frac{p_{\mu^\circ}(X^n) }{p_{\mu_1}(X^n)}\right)^{\frac{1}{2}}
\right] = e^{-\frac{n}{2} \renyi(\mu_1, \mu^{\circ})}.
\end{equation}
Since $M'_0$ and $M'_1$ are both CINECSI, it now follows from
Proposition~\ref{prop:withinconstants} that for some fixed $C_7 > 0$, 
\begin{equation}\label{eq:loting}
\renyi(\mu_1, \mu^{\circ})
\geq C_7 \| \mu_1 - \mu^{\circ} \|^2_2 \geq C_7 \| \mu_1 - \mu_0 \|^2_2,
\end{equation}
where the latter inequality follows by the definition of $\mu_0 =
\proj(\mu_1)$, see the explanation below \eqref{eq:mu0b}.  Combining
(\ref{eq:hoting}), (\ref{eq:joting}) and (\ref{eq:loting}), we have
thus shown that for all $n$, all $\mu_1 \in M_1$, all $A_2 \in \reals$,
\begin{equation}\label{eq:zoting}
\prob_{\mu_1}\left(\log{ \frac{p_{\mu_1}(X^n)
    }{p_{B',0}(X^n)}} < A_2  \right) \leq
C_6  e^{\frac{1}{2}A''}e^{-\frac{n}{2} C_7 \|\mu_1-  \mu_0\|_2^2}.
\end{equation}
The result now follows by combining (\ref{eq:start2}),
(\ref{eq:boting}) and (\ref{eq:zoting}).
\end{proof}

\section{Proof of Main Result, Theorem~\ref{thm:risk}}
\label{sec:mainproof}

\paragraph{Proof Idea}
The proof is based on analyzing what happens if $X_1, X_2, \ldots,
X_n$ are sampled from $p_{\mu^{(n)}_{1}}$, where $\mu^{(1)}_{1},
\mu^{(2)}_1, \ldots$ are a sequence of parameters in $M'_1$. We
consider three regimes, depending on how fast (if at all) $\mu_1^{(n)}$
converges to $\mu_0^{(n)}$ as $n \rightarrow \infty$. Here
$\mu_0^{(n)} = \proj(\mu_1^{(n)})$ is the projection of $\mu^{(1)}$
onto $M_0$, i.e. the distribution in $M_0$ defined, for each $n$,
as in (\ref{eq:mu0b}), with $\mu_1$ and $\mu_0$ in the role of
$\mu_1^{(n)}$ and $\mu_0^{(n)}$, respectively. Our regimes are defined
in terms of the function $f$  given by
\begin{equation}\label{eq:f}
f(n) := \frac{\sqd{\mbn}{\man}}{\frac{\log \log n}{n}} = 
 \frac{n \cdot \sqd{\mbn}{\man}}{\log \log n},
\end{equation}
which indicates how fast $\sq(\mbn , \man)$ grows relative to
the best possible rate $(\log \log n )/ n$. We fix appropriate
constants $\Gamma_1$ and $\Gamma_2$, and we distinguish, for all $n$ with $\Gamma_2 \log
n \geq  \Gamma_1$, the cases: 
$$
f(n)  \in \begin{cases} 
 [0, \Gamma_1] & \text{Case 1}  \\ [\Gamma_1,
 \Gamma_2 \log n] & \text{Case 2 (Theorem \ref{thm:upperbound})} \\ [\Gamma_2\log n, \infty] & \text{Case 3 (Theorem~\ref{thm:consistencycomplex})}.
\end{cases}
$$
For Case 1, the rate is easily seen to be upper bounded by $O((\log
\log n)/n)$, as shown inside the proof of Theorem~\ref{thm:risk}. In
Case 2, Theorem~\ref{thm:upperbound} establishes that the probability
that model ${\cal M}_0$ is chosen is at most of order $1/(\log n)$,
which, as shown inside the proof of Theorem~\ref{thm:risk}, again
implies an upper-bound on the rate-of convergence of $O((\log \log
n)/n)$.  Theorem~\ref{thm:consistencycomplex} shows that in Case 3,
which includes the case that $\sqd{\mbn}{\man}$ does not converge at
all, the probability that model ${\cal M}_0$ is chosen is at most of
order $1/n$, which, as again shown inside the proof of
Theorem~\ref{thm:risk}, again implies an upper-bound on the rate-of
convergence of $O((\log \log n)/n)$.

The two theorems take into account that $\mbn$ is not just a fixed
function of $n$, but may in reality be chosen by nature in a
worst-case manner, and that $f(n)$ may actually fluctuate between
regions for different $n$. Combining these two results, we finally
prove the main theorem, Theorem \ref{thm:risk}.

\begin{theorem}
\label{thm:consistencycomplex}
  Let $M_0, M_1$, $M'_1$ and $p_{\text{sw}, 1}(x^n)$ be as in Theorem~\ref{thm:risk}. Then there exist positive constants
  $\vec{b},\vec{c}$ such that for all $\mu_1 \in M'_1$, all $n$,
\begin{equation}\label{eq:expy}
\prob_{\mu_1}\left(\delta_\text{sw}(X^n) = 0 \right) \leq c_1 \cdot 
n^{m_1/2} \cdot e^{- c_2 n \| \mu_1^{(n)} - \mu_0^{(n)} \|_2^2} + \restje,
\end{equation}
where $\mu_0^{(n)} = \proj(\mu_1^{(n)})$ is as in (\ref{eq:mu0b}). As a
consequence, with $\Gamma_2  :=
c_2^{-1}(1+ m_1/2)$, we have the following: for every sequence
$\mu_1^{(1)}, \mu_1^{(2)}, \ldots$ with $f(n)$ as in (\ref{eq:f})
larger than $\Gamma_2 \log n$, we have 
$$
\prob_{\mu^{(n)}_1}\left(\delta_\text{sw}(X^n) = 0 \right) \leq  \frac{c_1}{n}
+ \restje.
$$
\end{theorem}

\begin{proof}
We can bound the probability of selecting the simple model by:
\begin{align*}
\prob_{\mu_1^{(n)}}\left(\delta_\text{sw}(X^n) = 0 \right) &= \prob_{\mu_1^{(n)}}\left(\switchratio \leq 1 \right)
= \prob_{\mu_1^{(n)}}\left(\frac{\sum_{i=0}^{ \infty }\pi(2^i)\bar p_{2^i}(X^n)  }{p_{B,0}(X^n)} \leq 1 \right)\\
&\leq \prob_{\mu_1^{(n)}}\left(\frac{ \pi(1) p_{B,1}(X^n)  }{p_{B,0}(X^n)} \leq 1 \right).
\end{align*}
Now (\ref{eq:expy}) follows directly by applying
Lemma~\ref{lem:bayesfactorrenyi} to the rightmost probability. 
For the second part, set $\Gamma_2  =
c_2^{-1}(1+ m_1/2)$. By assumption $f(n) > \Gamma_2 \log n$, we have
$\| \mu_1^{(n)}- \mu_0^{(n)} \|^2_2 > \Gamma_2 (\log n) (\log \log n)/n$. Applying  (\ref{eq:expy}) now gives the desired result.\end{proof}

\begin{theorem}\label{thm:upperbound}
  Let $f$ be as in (\ref{eq:f}) and $M'_1$ be as in Theorem
  \ref{thm:risk}. For any $\gamma > 0$, there exist constants $\Gamma_1, \Gamma_3 > 0$ such that, for
  every sequence $\mu_1^{(1)}, \mu_1^{(2)}, \ldots$ of elements of
  $M'_1$ with for all $n$, $f(n) > \Gamma_1$, we have
\begin{equation}\label{eq:start}
\prob_{\mu_1^{(n)}}\left( \switchratio \leq \gamma\right) \leq \frac{\Gamma_3}{\log n}.
\end{equation}
In particular, by taking $\gamma= 1$, we have
\begin{equation*}
{\mathbb P}_{\mu^{(n)}_1}(\delta_{\text{sw}}(X^n) = 0) \leq \frac{\Gamma_3}{\log n}.
\end{equation*}
The probabilities thus converge uniformly  at rate $O(1/(\log n))$ for all such sequences $\mu_1^{(1)}, \mu_1^{(2)}, \ldots$.
\end{theorem}

\begin{proof}
  We specify $\Gamma_1$ later. By assumption, we have
  $\pi(2^i) \gtrsim (\log n)^{-\kappa}$ for $i \in \{0, \ldots, \lfloor
  \log_2{n} \rfloor\}$.  We can restrict our attention to the strategy
  that switches to the complex model at the penultimate switching
  index, due to the following inequality: for any fixed $\gamma$, there exist positive constants $\vec{C}$  such that  for all large $n$:
\begin{align}
\prob_{\mu_1^{(n)}}\left( \switchratio \leq \gamma\right) &
\leq \prob_{\mu_1^{(n)}}\left( \frac{\sum_{i=0}^{\lfloor \log_2{n}\rfloor}\pi(2^i) \bar p_{2^i}(X^n) }{p_{B,0}(X^n)} \leq \gamma \right) \notag \\
&\leq  \prob_{\mu_1^{(n)}}\left( \frac{\sum_{i=0}^{\lfloor \log_2{n}\rfloor} \bar p_{2^i}(X^n) }{p_{B,0}(X^n)} \leq C_1 (\log{n})^\kappa  \right) \notag  \\
&\leq \prob_{\mu_1^{(n)}}\left( \frac{ \bar p_{2^{\lfloor \log_2{n}\rfloor - 1}}(X^n) } {p_{B,0}(X^n)} \leq  C_1 (\log{n})^\kappa  \right)\notag\\
&= \prob_{\mu_1^{(n)}}\left( \log{\frac{  \bar p_{2^{\lfloor \log_2{n}\rfloor - 1}}(X^n)} {p_{B,0}(X^n)}} \leq \kappa\log{  \log{n} } + C_2  \right) \label{eq:afschatting}. 
\end{align}
For the remainder of this proof, we will denote the penultimate switching index by $n^*$, that is: $n^* =  2^{\lfloor \log_2{n}\rfloor - 1}$. Now apply Lemma \ref{lem:8.1} twice, which gives that there exist $C_3,C_4$ such that, with probability at least $1- \restje$,
\begin{align}
 \log {\bar p_{{n^*}}(X^n)} &= 
\log p_{B,0}(X^{n^*}) + \log{p_{B,1}(X^n| X^{n^*})}=
\notag\\
&=
\log p_{B,0}(X^{n^*}) + \log{p_{B,1}(X^n)} - \log{p_{B,1}(X^{n^*})}\notag\\
&\geq \log p_{B,0}(X^{n^*}) + \log p_{\widehat{\mu}_1(X^n)}(X^n) -  \log p_{\widehat\mu_1(X^{n^*})}(X^{n^*})  + \frac{m_1}{2}\log{\frac{{n^*}}{n}}   - C_3 \notag\\ & \geq
\log p_{B,0}(X^{n^*}) + \log \frac{p_{\widehat{\mu}_1(X^n)}(X^n)}{p_{\widehat\mu_1(X^{n^*})}(X^{n^*})}   - C_4,
\end{align}
where we used that $\log \frac{n^*}{n}$ is of the order of a constant, because $n^*$ is between $\frac{n}{4}$ and $\frac{n}{2}$. From this, applying again Lemma~\ref{lem:8.1} twice, it follows that there exists $\vec{b}$ and $C_5, C_6$ such that for all $n$, with probability at least $1- \restje$,
\begin{align} 
\log \frac{\bar p_{{n^*}}(X^n)}{p_{B,0}(X^{n})}
&\geq  \log \frac{p_{B,0}(X^{n^*})}{p_{B,0}(X^{n})} + \log \frac{p_{\widehat{\mu}_1(X^n)}(X^n)}{p_{\widehat\mu_1(X^{n^*})}(X^{n^*})}   - C_4 \notag \\
& = 
- \log \frac{p_{\widehat{\mu}_0(X^n)}(X^n)}{p_{\widehat\mu_0(X^{n^*})}(X^{n^*})}  - \frac{m_0}{2}\log{\frac{{n^*}}{n}} + \log\frac{ p_{\widehat{\mu}_1(X^n)}(X^n)}{p_{\widehat\mu_1(X^{n^*})}(X^{n^*})}  - C_5 \notag \\
& \geq - \log \frac{p_{\widehat{\mu}_0(X^n)}(X^n)}{p_{\widehat\mu_0(X^{n^*})}(X^{n^*})} + \log\frac{ p_{\widehat{\mu}_1(X^n)}(X^n)}{p_{\widehat\mu_1(X^{n^*})}(X^{n^*})}  - C_6
\label{eq:ontwikkelingpn}
\end{align}
where we again used that $\log \frac{n^*}{n}$ can be bounded by constants. 
Let  $\cB_n$
be the event that (\ref{eq:ontwikkelingpn}) holds.  By
(\ref{eq:afschatting}) and (\ref{eq:ontwikkelingpn}), for all large
$n$, all $\beta \geq 1$,
\begin{align}\label{eq:basis}
& {\mathbb P}_{\mu_1^{(n)}}\left( \switchratio \leq \gamma 
\right)   \leq {\mathbb P}_{\mu_1^{(n)}}\left( \log{\frac{ \bar p_{{n^*}}(X^n)} {p_{B,0}(X^n)}} \leq \kappa\log{  \log{n} } + C_2  \right) \notag \\
\leq  & 
{\mathbb P}_{\mu_1^{(n)}}\left( \log{\frac{ \bar p_{{n^*}}(X^n)} {p_{B,0}(X^n)}} \leq \kappa\log{  \log{n} } + C_2 , \cB_n \right)+ 
{\mathbb P}_{\mu_1^{(n)}}\left( \cB^c_n \right)
\notag \\
 \leq & {\mathbb P}_{\mu_1^{(n)}}\left(- \log \frac{p_{\widehat{\mu}_0(X^n)}(X^n)}{p_{\widehat\mu_0(X^{n^*})}(X^{n^*})} + \log\frac{ p_{\widehat{\mu}_1(X^n)}(X^n)}{p_{\widehat\mu_1(X^{n^*})}(X^{n^*})}  - C_6 \leq \kappa\log{  \log{n} } + C_2 
\right) + \restje \notag \\
 = & {\mathbb P}_{\mu_1^{(n)}}\left( \cE^{(1)}_n  \right) + \restje
\leq  {\mathbb P}_{\mu_1^{(n)}}\left( \cE^{(\beta)}_n  \right) + \restje,
\end{align}
where we defined 
\begin{equation}\label{eq:cen}
  \cE^{(\beta)}_n =
  \left\{ \log
    \frac{p_{\widehat{\mu}_1(X^n)}(X^n)}{p_{\widehat\mu_1(X^{n^*})}(X^{n^*})} \cdot \frac{p_{\widehat\mu_0(X^{n^*})}(X^{n^*})}{p_{\widehat{\mu}_0(X^n)}(X^n)}  \leq  A^{(\beta)}_n \right\}
\end{equation}
and, for $\beta \geq 1$, we set $A^{(\beta)}_n = \beta \kappa \log
\log n + C_2 - C_6$.

Below, if a sample is split up into two parts $x_1, \ldots, x_{n^*}$
and $x_{n^*+1}, \ldots, x_n$, these partial samples will be referred
to as $x^{n^*}$ and $x^{>n^*}$ respectively.  We also suppress in our
notation the dependency of $A_n$, $\cE_n$ and $\cD_{j,n}$ as defined
below on $\beta$; all results below hold, with the same constants,
for any $\beta \geq 1$.

We will now bound the right-hand side of (\ref{eq:basis}) further.  Define the events
\begin{align*}
\cD_{1,n} & = \left\{ \log \frac{p_{\mu_1^{(n)}}(x^n)}{p_{\mu_1^{(n)}}(x^{n^*})}
\leq \log \frac{p_{\widehat{\mu}_1(X^n)}(x^n)}{p_{\widehat\mu_1(X^{n^*})}(x^{n^*})} + A_n
\right\} \\ 
\cD_{0,n} & = \left\{ \log \frac{p_{\mu_0^{(n)}}(x^n)}{p_{\mu_0^{(n)}}(x^{n^*})}
\geq \log \frac{p_{\widehat{\mu}_0(X^n)}(x^n)}{p_{\widehat\mu_0(X^{n^*})}(x^{n^*})} - A_n
\right\}.
\end{align*}
The probability in (\ref{eq:basis}) can be bounded, for all $\beta
\geq 1$, as 
\begin{align}\label{eq:decomp}
{\mathbb P}_{\mu_1^{(n)}}( \cE_n) & =  
{\mathbb P}_{\mu_1^{(n)}}( \cE_n, \cD_{0,n} \cap \cD_{1,n}) + {\mathbb P}_{\mu_1^{(n)}}( \cE_n, (\cD_{0,n} \cap \cD_{1,n})^c)  + \restje
\notag \\ & \leq  {\mathbb P}_{\mu_1^{(n)}}( \cE_n, \cD_{0,n}, \cD_{1,n}) 
+ {\mathbb P}_{\mu_1^{(n)}}(\cD_{1,n}^c) + {\mathbb P}_{\mu_1^{(n)}}(\cD_{0,n}^c) + \restje.
\end{align}
We first consider the first probability in (\ref{eq:decomp}): there are constants $\vec{C}$ such that, for all large $n$,  
\begin{align}\label{eq:bounda}
& {\mathbb P}_{\mu_1^{(n)}}( \cE_n, \cD_{0,n}, \cD_{1,n}) \notag \\ &\leq
{\mathbb P}_{\mu_1^{(n)}}\left(  
\log \frac{p_{{\mu}_1^{(n)}}(X^n)}{p_{\mu_1^{(n^*)}}(X^{n^*})} -  A_n 
+ \log \frac{p_{{\mu}_0^{(n)}}(X^{n})}{p_{{\mu}_0^{(n)}}(X^{n*})} -  A_n 
\leq  A_n \right) \notag \\  
&= {\mathbb P}_{\mu_1^{(n)}}\left( 
\log \frac{p_{\mu_1^{(n)}}(X^{> n^*})}{p_{\mu_0^{(n)}}(X^{> n^*})}  \leq 3 A_n \right)\notag \\
&\leq e^{ \frac{3}{2}    A_n} e^{- \frac{n}{4} \renyi( \mu_1^{(n)},
  \mu_0^{(n)})} \leq  e^{(3/2) \beta \kappa \log \log n + C_7} e^{- C_8 n \|  \mu_1^{(n)}-  \mu_0^{(n)} \|^2_2} 
=  e^{C_7} (\log n)^{(3/2) \beta \kappa - \Gamma_1\cdot C_8}, 
\end{align}
where $\Gamma_1$ is as in the statement of the theorem,  the second inequality follows by  Lemma \ref{lem:likbound}  and
noting $n^* < \frac{n}{2}$, we used Proposition~\ref{prop:withinconstants}.  

We now consider the second probability in (\ref{eq:decomp}). Using
$p_{\widehat{\mu}_1(X^n)}(x^n) \geq p_{\mu_1^{(n)}}(x^n)$ we have the
following, where we define the event $\cF_n = \{
\widehat{\mu}_1(X^{n^*}) \in M'_1 \}$ with $M'_1$ the CINECSI
subset
of $M_1$ mentioned in the theorem statement: there is $C_9, C_{10} > 0$
such that for al large $n$, 
\begin{align}\label{eq:boundb}
{\mathbb P}_{\mu_1^{(n)}}(\cD_{1,n}^c)
&=  {\mathbb P}_{\mu_1^{(n)}}\left(  \log \frac{p_{\mu_1^{(n)}}(X^n)}{p_{\mu_1^{(n)}}(X^{n^*})}
>  \log \frac{p_{\widehat{\mu}_1(X^n)}(X^n)}{p_{\widehat\mu_1(X^{n^*})}(X^{n^*})} +  A_n
\right)\notag\\ 
&\leq  
{\mathbb P}_{\mu_1^{(n)}}\left(  \log \frac{p_{\widehat{\mu}_1(X^n)}(X^n)}{p_{\mu_1^{(n)}}(X^{n^*})}
>  \log \frac{p_{\widehat{\mu}_1(X^n)}(X^n)}{p_{\widehat\mu_1(X^{n^*})}(X^{n^*})} +  A_n
\right) \notag\\
&\leq  {\mathbb P}_{\mu_1^{(n)}}\left(  
\log \frac{p_{\widehat\mu_1(X^{n^*})}(X^{n^*})}{p_{\mu_1^{(n)}}(X^{n^*})}
>    A_n  , \cF_n\right)  + {\mathbb P}_{\mu_1^{(n)}}\left(\cF_n^c
\right)
\notag \\
& \leq 
{\mathbb P}_{\mu_1^{(n)}}\left( D(\widehat{\mu}_1(X^{n^*}) \| \mu_1^{(n)} )>    A_n, \cF_n
\right) + \restje
\notag \\
& \leq 
{\mathbb P}_{\mu_1^{(n)}}\left( \|\widehat{\mu}_1(X^{n^*}) - \mu_1^{(n)} \|_2^2 >  C_9 A_n, \cF_n
\right) + \restje\notag \\
& \leq 
{\mathbb P}_{\mu_1^{(n)}}\left( \|\widehat{\mu}_1(X^{n^*}) - \mu_1^{(n)} \|_{\infty} >  \sqrt{C_9 A_n/m_1}
\right) + \restje\\ 
& \leq e^{- C_{10} A_n} = e^{- C_{10} (C_2 - C_6)} \frac{1}{(\log
  n)^{C_{10} \beta \kappa}},
\end{align}
where  we used  the KL robustness property 
(\ref{eq:robustness}), Proposition~\ref{prop:withinconstants} and Lemma~\ref{lem:klbound}. 

The third probability in  (\ref{eq:decomp}) is considered in a similar way. Using $p_{\widehat{\mu}_0(X^{n^*})}(X^{n^*}) \geq p_{\mu_0^{(n)}}(X^{n^*})$ we have $C_{11}, C_{12} > 0$ such that: 
\begin{align}\label{eq:boundc}
{\mathbb P}_{\mu_1^{(n)}}(\cD_{0,n}^c)
&=  {\mathbb P}_{\mu_1^{(n)}}\left(  \log \frac{p_{\mu_0^{(n)}}(X^n)}{p_{\mu_0^{(n)}}(X^{n^*})}
<  \log \frac{p_{\widehat{\mu}_0(X^n)}(X^n)}{p_{\widehat\mu_0(X^{n^*})}(X^{n^*})} - \frac{1}{3} A_n
\right)\notag\\ 
&\leq  
{\mathbb P}_{\mu_1^{(n)}}\left(  \log \frac{p_{\mu_0^{(n)}}(X^n)}{
p_{\widehat{\mu}_0(X^{n^*})}(x^{n^*})}
<  \log \frac{p_{\widehat{\mu}_0(X^n)}(X^n)}{p_{\widehat\mu_0(X^{n^*})}(X^{n^*})} - \frac{1}{3} A_n
\right) \notag\\
&= {\mathbb P}_{\mu_1^{(n)}}\left(  \log \frac{p_{\widehat{\mu}_0(X^n)}(X^n)}{p_{\mu_0^{(n)}}(X^n)} >   \frac{1}{3} A_n
\right)
\notag \\
& \leq C_{11} \frac{1}{(\log n)^{C_{12}\beta \kappa}} 
\end{align}
where we omitted the last few steps which 
are exactly as in (\ref{eq:boundb}).

We now finish the proof by combining (\ref{eq:decomp}),
(\ref{eq:bounda}), (\ref{eq:boundb}) and (\ref{eq:boundc}), which
gives that, if we choose $\beta \geq \max \{ 1/(\kappa C_{10}),
1/(\kappa C_{12}) \}$ and, for this choice $\beta$, we choose $\Gamma_1$ as in
(\ref{eq:bounda}) as $\Gamma_1 \geq (1 + (3/2) \beta \kappa)/C_8$, then we
have $\prob_{\mu_1^{(n)}}(\cE_n) \leq \Gamma_4 /(\log n)$ for some constant
$\Gamma_4$ independent of $n$; the result now follows from (\ref{eq:basis}).

\end{proof}

\noindent
\textbf{Proof of Theorem \ref{thm:risk}}
\begin{proof}
  We show the result in two stages. In Stage 1 we provide a tight
  upper  bound on the risk, based on an extension of the
  decomposition of the risk (\ref{eq:risk}) to general families and
  estimators $\breve\mu_0$ and $\breve\mu_1$ that are sufficiently
  efficient, i.e.  that satisfy (\ref{eq:suffeff}), and to losses
$\dgen{\cdot}{\cdot}$ equal to  squared error loss, standardized
  squared error loss and KL divergence (it is not sufficient to refer
  to Proposition \ref{prop:withinconstants} and prove the result only
  for squared error loss, because the equivalence result of
  Proposition \ref{prop:withinconstants} only holds on CINECSI sets
  and our estimators may take values outside of these; we do not need
  to consider R\'enyi and squared Hellinger divergences though,
  because these are uniformly upper bounded by KL divergence even for $\mu$
  outside any CINECSI set).  In Stage 2 we show how the bound implies the result. 

\commentout{
and 2 concern the
  upper-bound on the risk, which we will only show for
  $\dgen{\cdot}{\cdot}$ equal to squared error loss, standardized
  squared error loss and KL divergence. The upper bound for R\'enyi
  divergence and squared Hellinger divergence then follows
  automatically, since both are upper bounded by KL divergence. Stage
  3 concerns the lower-bound.}
  
\paragraph{Stage 1: Decomposition of Upper Bound on the Risk} 
Let $A_n$ be the event that $\cM_1$ is selected, as in Section
\ref{sec:Yang}. We will now show that, under the assumptions of Theorem~\ref{thm:risk}, we have for the constant $C$ appearing in
(\ref{eq:suffeff}), for all $\mu_1 \in M'_1$,
\begin{equation}\label{eq:trein}
R(\mu_1, \delta, n) \leq \frac{3C}{n} + 2 \prob(A^c_n) \dgen{\mu_1}{\mu_0},\end{equation}
where the left inequality holds for all divergence measures mentioned in the theorem, and the right inequality holds for $\dgen{\cdot}{\cdot}$ set to any of the squared error, the standardized squared error or the KL divergence. 

\commentout{{\em Proof of Stage 1, Left Inequality\/}
For any $\mu_1 \in M_1, \mu'\in M_0$, with $\mu_0 \in M_0$ as in (\ref{eq:mu0b}), we have
\begin{align}
R(\mu_1, \delta, n) &=  \E_{\mu_1}\left[ \1_{A_n} \dgen{\mu_1}{ \breve\mu_1(X^n)} + 
\1_{A^c_n} \dgen{\mu_1}{\breve\mu_0(X^n)}  \right] \nonumber \\
&\geq  \E_{\mu_1}\left[ \1_{A^c_n}  \left(\dgen{\mu_1}{\breve\mu_0(X^n)} \right)   \right]\notag\\
&\geq  \E_{\mu_1}\left[ \1_{A^c_n}  \left(\dgen{\mu_1}{\mu_0} \right)   \right]
 = \prob_{\mu_1}(A^c_n) \dgen{\mu_1}{\mu_0},
\end{align}
where the second inequality follows from the explanation below \eqref{eq:mu0b}.}

To prove (\ref{eq:trein}), we use that for the three divergences of interest,
for any $\mu_1 \in M_1, \mu \in M_0$, with $\mu_0 \in M_0$ as in (\ref{eq:mu0b}), we have
\begin{equation}\label{eq:dataday}
\dgen{\mu_1 }{\mu} \leq 2 ( \dgen{\mu_1}{ \mu_0} + \dgen{\mu_0 }{\mu}),
\end{equation}
For  $\dgen{\cdot}{\cdot}$ the KL divergence, this follows because 
\begin{align}
D(\mu_1 \| \mu) & = \E_{\mu_1}\left[- \log \frac{p_{\mu}(X)}{p_{\mu_1}(X)}\right] 
 = 
\E_{\mu_1}\left[- \log \frac{p_{\mu}(X)}{p_{\mu_0}(X)}\right] +  \E_{\mu_1}\left[- \log \frac{p_{\mu_0}(X)}{p_{\mu_1}(X)}\right] \notag \\
& = 
\E_{\mu_0}\left[- \log \frac{p_{\mu}(X)}{p_{\mu_0}(X)}\right] +  \E_{\mu_1}\left[- \log \frac{p_{\mu_0}(X)}{p_{\mu_1}(X)}\right],
\end{align}
where the last line follows by the robustness property of exponential
families (\ref{eq:robustnesspre}), since $\mu$ and $\mu_0$ are both in $M_0$. 

For $\dgen{\cdot}{\cdot}$ the  squared and standardized squared error case we show (\ref{eq:dataday}) as follows: 
Fix a matrix-valued function $J: M_1 \rightarrow \reals^{m_1^2}$ that maps each $\mu \in M_1$ to a positive definite matrix $J_{\mu}$. We can write
\begin{equation}\label{eq:dgen}
\dgen{\mu}{\mu'} = (\mu - \mu')^T J_{\mu} (\mu - \mu').
\end{equation}
where $J_{\mu}$ is the identity matrix for the squared error case, and
$J_{\mu}$ is the Fisher information matrix for the standardized
squared error case.  (\ref{eq:dataday}) follows since we can write,
for any function $J_{\mu}$ of the above type including these two:
\begin{align*}
(\mu_1 - \mu)^T J_{\mu_1} (\mu_1 - \mu) &= 
(\mu_1 - \mu_0 + \mu_0 -  \mu)^T J_{\mu_1} (\mu_1 - \mu_0 + \mu_0  -\mu) \\ &= 
(\mu_1 - \mu_0)^T J_{\mu_1} (\mu_1- \mu_0) +  (\mu_0 - \mu)^T J_{\mu_1} (\mu_0- \mu)
+ 2 (\mu_1 - \mu_0) J_{\mu_1} (\mu_0 - \mu) \\ & \leq 2 
\left((\mu_1 - \mu_0)^T J_{\mu_1} (\mu_1- \mu_0) +  (\mu_0 - \mu)^T J_{\mu_1} (\mu_0- \mu)
\right),
\end{align*}
where the last line follows because for general positive definite $m
\times m$ matrices $J$ and $m$-component column vectors $a$ and $b$,
$(b-a)^T  J (b-a) \geq 0$ so that $b^T J(b-a) \geq a^T J(b-a)$ and,
after rearranging, $b^T J b + a^T J a \geq 2 a^T J b$.

We have thus shown (\ref{eq:dataday}). It now follows that 
\begin{align}
 R(\mu_1, \delta, n) 
&=  \E_{\mu_1}\left[ \1_{A_n} \dgen{\mu_1}{ \breve\mu_1(X^n)} + 
\1_{A^c_n} \dgen{\mu_1}{\breve\mu_0(X^n)}  \right] \nonumber \\
&\leq  \E_{\mu_1}\left[ \dgen{\mu_1}{ \breve\mu_1(X^n)} + 
2 \cdot  \1_{A^c_n}  \left(\dgen{\mu_0}{\breve\mu_0(X^n)} + \dgen{\mu_1}{\mu_0} \right)   \right]
 \nonumber \\
& \leq \frac{3C}{n} + 2 \prob(A^c_n) \dgen{\mu_1}{\mu_0},  \label{eq:extendedrisk}
\end{align}
where we used (\ref{eq:dataday}) and our condition (\ref{eq:suffeff}) on $\breve\mu_0$ and $\breve\mu_1$. We have thus shown (\ref{eq:trein}).

\paragraph{Stage 2}
\commentout{
{\em Proof of Lower-Bound\/} From the first inequality of (\ref{eq:trein}) and the fact that $M'_1$ is CINECSI, there exists a constant $C_2$ such that, with $R$ defined relative to the divergence $\dgen{\cdot}{\cdot}$ of interest, for every sequence $\mu_1^{(1)}, \mu_1^{(2)}, \ldots$ with each $\mu_1^{(i)} \in M'_1$, 
$$
R(\mu^{(n)}_1, \delta, n)\geq  C_2 \cdot
\prob_{\mu^{(n)}_1}(A^c_n) \sqd{\mu^{(n)}_1}{\mu^{(n)}_0}.$$ 
The lower bound follows by taking a specific sequence $\mu_1^{(1)}, \mu_1^{(2)}, \ldots$ such that with $f(n)$ as in (\ref{eq:f}) and $c_1$ as in Theorem~\ref{lem:boundswitchratio}, $f(n) = c_1/2$.  \stephanie{Waarom $c_1/2$? Bedoel je misschien $c_3/2$?}Applying 
the second part of that Theorem it follows that there is a $c_2 > 0$ such that for all large $n$,  $R(\mu^{(n)}_1, \delta, n) \geq (c_2/2) \cdot (\log \log n)/n$, \stephanie{Waarom hier dan delen door twee? De constante wordt $C_2\cdot(c_3/2)$?}which proves the lower-bound.

{\em Proof of Upper-Bound\/}} We proceed to prove our risk upper bound
for the squared error loss, standardized squared error loss and KL
divergence, for which the right inequality in (\ref{eq:trein}) holds;
the result then follows for squared Hellinger and R\'enyi divergence
because these are upper bounded by KL divergence. From
(\ref{eq:trein}) we see that it is sufficient to show that for
all $n$ larger than some $n_0$,
\begin{equation}\label{eq:cs}
\sup_{\mu_1 \in M'_1} \{ \prob_{\mu_1}(A^c_n) \dgen{\mu_1}{\mu_0} \} = O\left(\frac{\log \log n}{n} \right),
\end{equation}
for our three choices of $\dgen{\cdot}{\cdot}$. 
We first note that, since $M'_1$ is CINECSI, $\sup_{\mu_1 \in M'_1}
\dgen{\mu_1}{\mu_0}$ is bounded by some constant $C_1$. It thus
follows by Proposition~\ref{prop:niceml} that there exists some CINECSI subset $M''_1$
of $M_1$ such that, with $B^c_n \subset A^c_n$ defined as $B_n^c =
\{x^n: \delta(x^n) = 0 ; \widehat\mu_1(X^n) \in M''_1\}$, we have
\begin{align*}
\sup_{\mu_1 \in M'_1} \{ \prob_{\mu_1}(A^c_n) \dgen{\mu_1}{\mu_0} \} &=
\sup_{\mu_1 \in M'_1} \{ (\prob_{\mu_1}(B^c_n) + \prob_{\mu_1}(A^c_n \setminus 
B^c_n) ) \dgen{\mu_1}{\mu_0} \} \\ &=
\sup_{\mu_1 \in M'_1} \{ \prob_{\mu_1}(B^c_n) \dgen{\mu_1}{\mu_0} \} + 
C_1 \cdot \prob_{\mu_1}(\widehat{\mu}^{(1)} \not \in M''_1) \\ &= 
\sup_{\mu_1 \in M'_1} \{ \prob_{\mu_1}(B^c_n) \dgen{\mu_1}{\mu_0} \}+ 
\restje,
\end{align*}
so that it is sufficient if we can show (\ref{eq:cs}) with $B^c_n$ instead of $A^c_n$. But on the set $B^c_n$, all three divergence measures considered are within constant factors of each other, so that it is sufficient if we 
can show that there is a constant $C_2$ such that for all $n$ larger than some $n_0$, 
\begin{equation}\label{eq:ds}
\sup_{\mu_1 \in M'_1} \{ \prob_{\mu_1}(B^c_n) \cdot \sqd{\mu_1}{\mu_0} \} \leq C_2 \cdot \frac{\log \log n}{n}. 
\end{equation}
Now, fix some $\mu_1 \equiv \mu_1^{(n)}$ and consider $f(n)$ as in
(\ref{eq:f}).  By Theorem~\ref{thm:consistencycomplex},
$\prob_{\mu_1}(B^c_n) \leq C_3/n$ for some constant $C_3$ that can be
chosen uniformly for all $\mu_1 \in M'_1$ whenever $f(n) > \Gamma_2 \log n$
with $\Gamma_2$ as in that theorem. Using also that $\| \mu_1 - \mu_0
\|^2_2$ is bounded by $C_1$ as above, it follows that (\ref{eq:ds})
holds whenever $f(n) > \Gamma_2 \log n$ and $ (C_1 C_3) /n \leq C_2 (\log
\log n)/n$, i.e. whenever $f(n) > \Gamma_2 \log n$ and $C_2 \geq C_1 C_3 /
(\log \log n)$. 

Second, suppose that $\Gamma_1 < f(n) \leq  \Gamma_2 \log n$ with $\Gamma_1$ as in 
Theorem~\ref{thm:upperbound}. Then by that theorem, uniformly 
for all $\mu^{(n)}_1$ with such $f(n)$, we have, with $\Gamma_3$ as in that theorem,
\begin{align*}
& \sqd{\mu^{(n)}_1}{\mu^{(n)}_0} \cdot {\mathbb P}_{\mu^{(n)}_1}(\delta_{\text{sw}}(X^n) = 0)
= f(n) \cdot \frac{\log \log n}{n} \cdot {\mathbb
  P}_{\mu^{(n)}_1}(\delta_{\text{sw}}(X^n) = 0) 
\leq  \\ &  \Gamma_2 \cdot (\log n)  \cdot  \frac{\log \log n}{ n} \cdot
{\mathbb P}_{\mu^{(n)}_1}(\delta_{\text{sw}}(X^n) = 0) \leq \Gamma_2\Gamma_3 
\cdot  \frac{\log \log n}{ n},
\end{align*}
where $\mu^{(n)}_0 = \proj(\mu^{(n)}_1)$ is defined as in
(\ref{eq:mu0b}), so that (\ref{eq:ds}) holds again whenever  $C_2 \geq
\Gamma_2\Gamma_3$. 

Finally, suppose that $f(n) \leq\Gamma_1$ with $\Gamma_1$ as in 
Theorem~\ref{thm:upperbound}. Then (\ref{eq:ds}) holds whenever $C_2
\geq \Gamma_1$. Combining the three cases we find that (\ref{eq:ds}) holds
whenever $C_3 \geq \max \{\Gamma_1, \Gamma_2\Gamma_3, \allowbreak C_1C_3/(\log \log n) \}$; the
result is proved.

\end{proof}

\section{Switching as in Van Erven et al. (\citeyear{Erven2012})}
\label{app:realswitch}

The basic building block of the switch distribution and criterion as
formulated by Van Erven et al. (\citeyear{Erven2012}) is a  countable set of \emph{sequential
  prediction strategies\/} (also known as `prequential forecasting systems' \citep{Dawid84})
$\{p_k \mid k \in \mathcal{K}\}$, where $\mathcal{K}$ is a finite or
countable set indexing the basic models under consideration. Thus,
each model is associated with a corresponding prediction strategy,
where a prediction strategy $p$ is a function from $\bigcup_{i \geq 0}
\mathcal{X}^i$ to the set of densities on $\mathcal{X}$, where
$p(\cdot \mid x^{n-1})$ denotes the density on $\mathcal{X}$ that
$x^{n-1}$ maps to, and $p(x_n \mid x^{n-1})$ is to be interpreted as
the probabilistic prediction that strategy $p$ makes for outcome $X_n$
upon observation of the first $n-1$ outcomes, $X^{n-1} = x^{n-1}$.  For
example, for a parametric model $\{p_\theta \mid \theta \in \Theta\}$ one
can base $p_k$ on a Bayesian marginal likelihood, $p_B(x^n) :=
\int_\Theta \omega(\theta) p_\theta(x^n)d\theta$, where $\omega$ is a
prior density on $\Theta$. The corresponding prediction strategy could
then be defined by setting $p_k(x_{n} \mid x^{n-1}) :=
p_B(x^n)/p_B(x^{n-1})$, the standard Bayesian predictive
distribution. In this paper, the basic strategies $p_k$ were always
Bayesian predictive distributions, but, in the spirit of
\cite{Dawid84}, one may consider other choices as well.

After constructing the set of basic prediction strategies, a new
family of prediction strategies that switch between the strategies in
the set $\{p_k \mid k \in \mathcal{K}\}$ is defined. Formally, let
$\mathbb{S}$ be the set
\begin{equation}\label{eq:s}
\mathbb{S} = \left\{ \left((t_1, k_1), \ldots, (t_m, k_m) \right) \in \left(\mathbb{N} \times \mathcal{K}\right)^m | m \in \mathbb{N}, 1 = t_1 < t_2 < \ldots < t_m \right\}.
\end{equation}
Each $s \in \mathbb{S}$ specifies the times $t_1, \ldots, t_m$ at which a switch is made between the prediction strategies from the original set, identified by the indices $k_1, \ldots, k_m$. The new family $Q = \{q_s \mid s \in \mathbb{S}\}$ is then defined by setting, for all $n, x^n \in {\cal X}^n$: 
\begin{equation}\label{eq:seqnature}
q_s(x_{n} \mid x^{n-1}) = p_{k_j}(x_n \mid x^{n-1}), \quad t_j \leq n  < t_{j+1},
\end{equation}
with $t_{m+1} = \infty$ by convention. We now define $q_s(x^n) =
\prod_{i=1}^n q_s(x_i \mid x^{i-1})$; one easily verifies that this
defines a joint probability density on ${\cal X}^n$.  

We now place a prior mass function $\pi'$ on $\mathbb{S}$ and define,
for each $n$, the {\em switch distribution\/} in terms of its joint
density for ${\cal X}^n$ and $\mathbb{S}$:
\begin{equation*}
  p_{\text{sw}}(x^n,s) = q_s(x^n)\pi'(s), \quad  p_{\text{sw}}(x^n) = \sum_{s \in \mathbb{S}} p_{\text{sw}}(x^n,s) = \sum_{s \in \mathbb{S}} q_s(x^n)\pi'(s).
\end{equation*}
If the $p_k$ are defined as Bayesian predictive distributions as
above, then, as explained by Van Erven et al. (\citeyear{Erven2012}), the density
$p_{\text{sw}}(x^n)$ can be interpreted as a Bayesian marginal density
of $x^n$ under the prior $\pi'$ on meta-models (model sequences) in
$\mathbb{S}$.

The switch distribution can be used to define a model selection criterion $\delta'_{\text{sw}}$ by selecting the model with highest posterior probability under the switch distribution. This is done by defining the random variable $K_{n+1}(s)$ on $\mathbb{S}$ to be the index of the prediction strategy that is used by $q_s$ to predict the $(n+1)$th outcome. The model selection criterion is then:
\begin{align}
\delta'_\text{sw}(x^n) &= \arg\max_k p_\text{sw}(K_{n+1} = k \mid x^n) = \arg\max_k \frac{\sum_{s: K_{n+1}(s) = k}p_\text{sw}(x^n, s)}{p_\text{sw}(x^n)}\notag\\
&= \arg\max_k \frac{ \sum_{s: K_{n+1}(s) = k}q_s(x^n)\pi'(s)}{\sum_{s \in \mathbb{S}}q_s(x^n)\pi'(s)}, \label{eq:deltaswold}
\end{align}
with ties resolved in any way desired.

In our nested two-model case, one might use, for example, a prior
$\pi'$ with support on 
$${\mathbb S}' = \{(1,0),(1,1),(\; (1,0),(2,1)\;), (\; (1,0),(4,1)\;),
(\; (1,0),(8,1)\;), (\; (1,0),(16,1), \ldots\;) \}.$$ Such a prior
expresses that at time $1$, for the first prediction, one can either
switch to (i.e., start with), model $0$, and keep predicting according
to its Bayes predictive distribution --- this strategy gets weight
$\pi((1,0))$. Or one can start with model $1$, and keep predicting
according to its Bayes predictive distribution --- this strategy gets weight
$\pi((1,1))$. Or one can start with model $0$ and switch to model $1$
after $2^i$ observations and then stick with $1$ forever --- this
strategy gets weight $\pi((\; (1,0), (2^i,1) \;))$. If we now start
with a prior $\pi$ on $\{1,2, \ldots\}$ as in the main text and define
$\pi'((1,0)) = 1/2$, $\pi'((1,1)) = (1/2) \cdot \pi(1)$, and for $i \geq
1$, $\pi'((1,0),(2^i,1)) = (1/2) \cdot  \pi(2^i)$, then $\sum_{s \in {\mathbb
    S}'} \pi'(s) = 1$, so $\pi'$ is a probability mass function.  A
simple calculation gives that (\ref{eq:deltaswold}) based on switch
prior $\pi'$ now chooses model $1$ if
\begin{equation}\label{eq:vroeger}
{\sum_{1 \leq t < n} \bar p_t(x^n)\pi(t)} > {(1 +g(n) ) \cdot  p_{B,0}(x^n)},
\end{equation}
where $g(n) = \sum_{t \geq n}
  \pi(t)$; note that $g(n)$ is decreasing and converges to $0$ with
  increasing $n$.  
(\ref{eq:vroeger}) is thus an instance of the switch criterion of Van Erven et al. (\citeyear{Erven2012}).
Comparing this to (\ref{eq:deltasworig}), the criterion used in this
paper, after rearranging we see that it chooses model $1$ if
$$
{\sum_{1 \leq t < n} \bar p_t(x^n)\pi(t)} > (1- g(n)) \cdot  { p_{B,0}(x^n)},
$$
which is more likely by constant factor to select model $\M_0$, the
factor however tending to $1$ with increasing $n$. It is completely
straightforward to check that Theorem~\ref{thm:risk} and all other
results in this paper still hold if $\delta_\text{sw}$ with prior
$\pi$ as in the main text is replaced by $\delta'_\text{sw}$ with
corresponding prior $\pi'$ as defined here; thus our results carry
over to the original definitions of Van Erven et al.
(\citeyear{Erven2012}). Similarly, the proof for the strong
consistency of $\delta'_\text{sw}$ given by Van Erven et
al. (\citeyear{Erven2012}) carries through for $\delta_{\text{sw}}$,
needing only trivial modifications. From (\ref{eq:vroeger}) we see
that modifying the prior $\pi$ in either our or Van Erven et al.'s
original criterion has a similar effect as keeping the same $\pi$ but
switching between the two versions of the switch criterion.

\bibliographystyle{plainnat}
\bibliography{switch}
\end{document}